\documentclass[12pt]{iopart}
\usepackage{iopams}  
\usepackage{amsopn}
\usepackage{hyperref}
\usepackage{algorithm2e}
\usepackage{epsfig}
\usepackage{graphicx}

\newcommand{\B}[1]{{\bf #1}}
\newcommand{\Sc}[1]{{\mathcal{#1}}}

\newtheorem{lemma}{Lemma}[section]
\newtheorem{remark}[lemma]{Remark}
\newtheorem{theorem}[lemma]{Theorem}
\newtheorem{corollary}[lemma]{Corollary}

\def\realbold#1{{\mbox{\boldmath $\bf #1$}}}

\def\@undertilde#1{\oalign{{\realbold{#1}}\crcr\hidewidth
\vbox to .2ex{\hbox{\bf\char126}\vss}\hidewidth}}

\DeclareMathOperator*{\argmin}{argmin}

\begin{document}

\title{Estimating Nuisance Parameters in Inverse Problems}

\author{Aleksandr Y. Aravkin and Tristan van Leeuwen}

\address{Department of Earth, Ocean and Atmospheric Sciences,
The University of British Columbia,
6339 Stores Road,
Vancouver, British Columbia
Canada, V6T 1Z4.}
\ead{\{saravkin, tleeuwen\}@eos.ubc.ca}

\begin{abstract}
Many inverse problems include nuisance parameters which, while not of direct interest, 
are required to recover primary parameters. 
Structure present in these problems allows efficient optimization strategies --- a well known example 
is {\it variable projection}, where nonlinear least squares problems which are linear in some parameters
can be very efficiently optimized. 
In this paper, we extend the idea of {\it projecting} out a subset over the variables 
to a broad class of maximum likelihood (ML) and maximum {\it a posteriori} likelihood
 (MAP) problems with nuisance parameters, such as variance or degrees of freedom. 
 As a result, we are able to incorporate nuisance parameter estimation into 
large-scale constrained and unconstrained inverse problem formulations. 
We apply the approach to a variety of problems, 
including  estimation of unknown variance parameters in the Gaussian model,  
degree of freedom (d.o.f.) parameter estimation in the context of robust inverse problems, automatic calibration, 
and optimal experimental design. 
Using numerical examples, we demonstrate improvement in recovery of primary parameters 
for several large-scale inverse problems.
The proposed approach is compatible with a wide variety of algorithms and formulations, and 
its implementation requires only minor modifications to existing algorithms.   
\end{abstract}

\maketitle

\section{Introduction}

Many inverse problems can be formulated as optimization problems 
of the form 
\begin{equation}
\label{InverseProblemClass}
\Sc{P} \quad \min_{x \in \Sc{X}, \theta} g(x,\theta)\;,
\end{equation}
where $g: \mathbb{R}^n \times \mathbb{R}^k \rightarrow \mathbb{R}$ is
a twice differentiable function, $\Sc{X} \subset \mathbb{R}^n$, $x$ is a primary set of 
parameters of interest, while $\theta \in \mathbb{R}^k$ is a secondary
set of nuisance parameters, such as variance parameters, application-specific 
tuning parameters, regularization parameters, or degrees of freedom 
parameters. In many settings, $k\ll n$.

A rich source of examples in~\eref{InverseProblemClass} 
is the class of separable least-squares problems,  extensively studied  
over the last 40 years~\cite{GolubPereyra,GolubPereyra2003,Osborne2007}.
A problem in this class is given by 
\begin{equation}
\label{SeparableLeastSquares}
\min_{x, \theta} \|y - \Phi(x)\theta\|_2^2\;, 
\end{equation}
where the matrix $\Phi(x)$ is parametrized by $x$. Note that $g$ has a very 
special form in this case, and $\Sc{X} = \mathbb{R}^n$.
For problems in this class, the major
insight is to exploit the structure of the problem to obtain a reduced
problem
\begin{equation}
\label{modProblem}
\min_{x} \|y - \Phi(x)\bar\theta(x)\|_2^2\;, 
\end{equation} 
where
\begin{equation}
\label{thetaOfx}
\bar\theta(x) = \argmin_\theta \|y - \Phi(x)\theta\|_2^2\;.
\end{equation}
At first glance, this does not make the problem easier to solve. However, it
turns out that~\eref{modProblem} can be solved using black-box approaches as long as we
re-evaluate $\bar\theta(x)$ for any given $x$, but treat $\theta$ as fixed 
whenever $x$ is updated. The problem~\eref{thetaOfx} has a closed form solution,
and as noted in~\cite{GolubPereyra2003}, this approach converges much faster
than optimization approaches to minimize the full functional~\eref{SeparableLeastSquares}
using descent methods for $(x, \theta)$.

In this paper, we consider problems of type~\eref{InverseProblemClass},
where we can easily compute $\bar\theta(x) = \argmin_\theta g(x,\theta)$.
We show that many algorithms for solving instances of~\eref{InverseProblemClass} with $\theta$ fixed
can be easily modified to solve the joint inverse problem in $x$ and $\theta$.   
We provide explicit details for several important classes of problems in~\eref{InverseProblemClass}, 
including variance and degrees of freedom (d.o.f.) estimation, and automatic calibration 
of nonlinear least squares and robust inverse problem formulations.  

The paper proceeds as follows. In Section \ref{SecGeneral}, we review the necessary
theory underlying our approach to the entire class~\eref{InverseProblemClass}.
In Section~\eref{MapEst}, we discuss the role of nuisance parameters, such 
as variance and degrees of freedom, in MAP estimation formulations. 
We present two important applications in detail:
\begin{enumerate}
\item Variance estimation in multiple data sets (see~\cite{BellBurkeSchumitzky}). 

\item Estimation of variance and d.o.f. for Student's t formulations (see~\cite{Lange1989}). 
\end{enumerate} 
Both are illustrated on a seismic imaging problem where the data are contaminated 
with various types of noise.

In Section \ref{RobustSourceParam}, we discuss the automatic calibration problem, 
where the forward model includes a calibration parameter that needs to be estimated.
We illustrate the approach on a seismic imaging problem where the calibration parameter
are frequency-dependent source-weights.
We discuss the application of the proposed approach to Optimal Experimental Design in section
\ref{OED}.

Finally, we discuss other possible applications and present conclusions.


\section{General Formulation}
\label{SecGeneral}

We consider problems of the form~\eref{InverseProblemClass},
and assume that for any given $x \in \Sc{X}$, one can easily find 
\begin{equation}
\label{Projection}
\bar\theta(x) \in \argmin_{\theta} g (x, \theta)\;.
\end{equation}
This condition can be relaxed, and $\bar\theta(x)$ can be considered a local minimum.
Rather than working to solve~\eref{InverseProblemClass}, we can instead focus on 
the {\it reduced objective} 
\begin{equation}
\label{ProjectedClass}
\tilde g(x) = g(x, \bar\theta(x))\;. 
\end{equation}
This approach is justified by the following theorem, adapted from~\cite[Theorem 2]{Bell2008ADo}.

\begin{theorem}
\label{ProjectionTheorem}

Suppose that $\Sc{U}\subset\mathbb{R}^n$ and 
$\Sc{V} \subset \mathbb{R}^k$ are open, and
$g(x,\theta)$ is twice continuously differentiable on $\Sc{U}\times\Sc{V}$. 
Define the optimal value function 
\begin{equation}
\label{valueFunction}
\tilde g(x) = \min_\theta g(x, \theta)\;.
\end{equation}
Suppose that $\bar x \in \Sc{U}$ and $\bar \theta \in \Sc{V}$ are such 
that $\nabla_\theta g(\bar x, \bar \theta) = 0$ and $\nabla_\theta^2 g(\bar x, \bar \theta)$
is positive definite. Then there exist neighbourhoods of $\bar x$ and $\bar \theta$
and a twice continuously differentiable function $\bar\theta: \Sc{U} \rightarrow \Sc{V}$
where $\bar\theta(x)$ is the unique minimizer of $g(x, \cdot)$ on $\Sc{V}$. 

Then $\tilde g(x)$ is twice continuously differentiable, with 
\begin{eqnarray}
\label{firstProjDer}
\nabla_x \tilde g(\bar x) &=& \nabla_x g(\bar x, \bar\theta(\bar x))\\
\label{secondProjDer}
\nabla_x^2 \tilde g(\bar x) &=& \nabla_x^2 g(\bar x, \bar\theta(\bar x)) +
\nabla^2_{x,\theta} g(\bar x, \bar\theta(\bar x))\nabla_x\bar\theta(\bar x)\;.
\end{eqnarray}
\end{theorem}

\begin{remark}
\label{PonyRemark}
Theorem~\ref{firstProjDer} provides {\emph sufficient} conditions for existence of the first and second derivatives 
of $\tilde g$. In practice, these derivatives may exist even if the smoothness hypotheses are not satisfied. 
Consider $g(\theta, x) = \frac{x^4}{2} + \theta^2 - |\theta|x^2$. In this case, $|\bar \theta(x)| = \frac{x^2}{2}$, 
so $\tilde g(x) = \frac{x^4}{4}$ is smooth even though $g(x,\theta)$ is not.  
\end{remark}

Theorem~\ref{firstProjDer} suggests a natural approach to designing algorithms for minimizing
$\tilde g(x)$. 
In the unconstrained case (i.e. $\Sc{X}$ is the whole space), 
consider iterative methods of the form 
\begin{equation}
\label{IteratedMethod}
x^{k+1} = x^k - \gamma_k H_k^{-1} \nabla_x \tilde g(x^k) = x^k - \gamma_k
H_k^{-1} \nabla_x g(x^k, \bar\theta(x^k))\;.
\end{equation}
Specifically, $H_k = I$ yields Cauchy's steepest descent,
$H_k = \nabla_x^2 \tilde g(x^k)$ yields a modified 
Newton method, and approximations to $H_k$ that use 
only first order derivative information yield Gauss-Newton or 
Levenberg-Marquardt type methods. 
A quasi-Newton method such as BFGS or L-BFGS may be similarly implemented
using only information from~\eref{firstProjDer}.
  
If $\Sc{X}$ is a closed and bounded set that allows a simple projection, 
such as a set of box constraints ($\{x: l \leq x \leq u\}$, 
an ellipsoidal set $\{x: \|x\|_M \leq \tau\}$, or the $1$-norm ball ($\{x: \|x\|_1 \leq \tau\}$), 
this can be exploited to solve~\eref{InverseProblemClass}. 
For example, we can use a modified projected gradient method
\begin{equation}
\label{ProjectedGradient}
x^{k+1} = \Sc{P}_{\Sc{X}}[x^k - \gamma  \nabla_x \tilde g(x^k)]
\end{equation}
or an appropriately modified projected quasi-Newton method, such as the one 
described in~\cite{Schmidt09optimizingcostly}. The point is that the structure of 
$\Sc{X}$ does not enter into the computation of~\eref{firstProjDer} or~\eref{secondProjDer}, 
so a natural strategy is to compute these quantities first and then 
apply methods that exploit the structure of $\Sc{X}$. Moreover, 
we show in the next corollary that the point $(\bar x, \theta(\bar x))$  
satisfies the first order necessary conditions for the original (constrained) problem. 

\begin{corollary}
\label{ConstrainedCorollary}
Suppose the hypotheses of Theorem~\ref{firstProjDer} hold, and the additional
constraint $x \in \Sc{X}$ is imposed, where $\Sc{X}$ is a closed convex set. 
If $\bar x$ satisfies the first order necessary  conditions for $\tilde g(x)$, 
then $(\bar x, \bar\theta)$ with $\bar \theta = \theta(\bar x)$ satisfies
the first order necessary conditions for $g(x, \theta)$. 
\end{corollary} 

\paragraph{Proof:}
The first order necessary conditions for~\eref{InverseProblemClass}
are 
\begin{equation}
\label{NecessaryConditions}
\begin{array}{ccc}
\nabla_\theta g(\bar x,\bar \theta) &=& 0\\
-\nabla_x g(\bar x,\bar\theta) &\in& N_{\Sc{X}}(\bar x)
\end{array}
\end{equation}
where $N_{\Sc{X}}(\bar x)$ is the normal cone to $\Sc{X}$ at the point $\bar x$ (see~\cite{RTR} for details). 
The first order necessary condition for $\bar x$ to be a minimizer of the reduced objective~\eref{ProjectedClass}
is  
\begin{equation}
\label{NecessaryConditionsP}
-\nabla_x \tilde g(\bar x) \in N_{\Sc{X}}(\bar x)
\end{equation}
Since we have $\nabla_x \tilde g(\bar x) =\nabla_x g(\bar x,\bar\theta) $ by Theorem~\ref{firstProjDer}, 
$(\bar x, \theta(\bar x))$ satisfies~\eref{NecessaryConditions}  
if and only if $\bar x$ satisfies~\eref{NecessaryConditionsP}. 
On the other hand, $\theta(\bar x)$ satisfies the first equation of~\eref{NecessaryConditions} by 
construction.


Thus, for many applications (both constrained and unconstrained), 
we can systematically extend many standard algorithms  
for minimizing $g(x, \theta)$ with $\theta$ fixed
to extended problems~\eref{InverseProblemClass}. 
This approach avoids computing the full Hessian of the modified objective~\eref{secondProjDer}, 
since it involves $\nabla_x\bar\theta(x)$.

In the next sections, we present some of these applications
and provide full algorithmic details and numerical work. 

\section{Complicating Parameters in Maximum Likelihood Estimation}
\label{MapEst}

Many inverse problems can be formulated as maximum likelihood (ML) problems
within a statistical modeling framework. Given data $y$, we want to solve for 
parameters of interest $x$, using the fact that the parameters are related to the data 
via a (possibly nonlinear) forward model:
\begin{equation}
\label{MLformulation}
d = F(x) + \epsilon\;.
\end{equation}
The $\epsilon$ term in~\eref{MLformulation} reflects a statistical model of the discrepancy
between the model $F(x)$ and the true data $d$. Independent, identically distributed (i.i.d.) Gaussian 
errors $\epsilon \sim \B{N}(0, \sigma^2I)$ are a common choice, 
and even though the variance parameter $\sigma^2$ is unknown, 
it does not affect the maximum likelihood formulation in $x$. 
This is not true if the data come from different sources, with each group having its own 
parameter $\sigma^2_i$. 

More generally, $\epsilon_i$ may come from a range of parametric distributions. 
The Student's t distribution has been applied in many instances where large measurement errors 
are common or unexplained artifacts in the data are an issue~\cite{aravkin12,AravkinFriedlanderHerrmannVanleeuwen:2011,Lange1989}. 
These applications require estimates for degrees of freedom and variance parameters 
even with the i.i.d. assumption on the errors.

If we take $\theta$ to be unknown nuisance parameters, the general maximum likelihood formulation
for estimating $x$ in model~\eref{MLformulation} takes the form~\eref{InverseProblemClass}. 
The method proposed in this paper is well suited for online estimation of $\theta$, 
and in the remainder of the section we provide full exposition for the multiple sources of error 
example and for Student's t hyperparameter estimation.

\subsection{Variances in Multiple Data Sets}
\label{VarEst}

Estimating variance parameters in multiple datasets is an important problem in 
many areas, including drug and tracer kinetics~\cite{BellBurkeSchumitzky}, and geophysics. 
In this section, we review the formulation presented in~\cite{BellBurkeSchumitzky}, 
and show that the algorithm derived in~\cite{BellBurkeSchumitzky} 
follows immediately from the general approach we propose here, i.e. 
it is a Gauss-Newton method of form~\eref{IteratedMethod}. 
We present a numerical example, illustrating the importance of variance parameter estimation 
for a large-scale geophysical inverse problem. 
We also extend the approach to the (fully observed) multivariate Gaussian case with 
correlations between measurement errors.

We are given $M$ experiments indexed by $i$, each of which yields $N_i$ measurements 
and has its own variance parameter $\sigma_i$. All experiments share a common set of primary parameters $x$: 
\begin{equation}
\label{statisticalModel}
d_{i} = F_{i}(x) + \epsilon_{i}
\end{equation} 
where $d_i \in \mathbb{R}^{N_i}$, $F_i(x)$ is the modeling operator for the
$i^{\mathrm{th}}$ experiment and $\epsilon_{i} \sim \B{N}(0, \sigma^2_iI)$. 
If the variance parameters are fixed, the ML estimation problem for $x$
is given by
\begin{equation}
\label{NLLS}
\min_{x} \quad
\sum_{i=1}^M 
\frac{1}{\sigma^2_i}\|d_i - F_i(x) \|_2^2\;.
\end{equation}
The joint ML estimation problem for $x$ and $\sigma^2 = \{\sigma_i^2\}$ is given by 
\begin{equation}
\label{MAPmultiple}
\min_{\sigma^2, x} \quad
g(x, \sigma^2) :=  
\sum_{i=1}^M \left( 
N_i \log(2\pi \sigma^2_i) 
+ \frac{1}{\sigma^2_i}\|d_i - F_i(x) \|_2^2
\right)\;.
\end{equation}
This is a special example of~\eref{InverseProblemClass}. 

With $x$ fixed,~\eref{MAPmultiple} separates, 
and~\eref{Projection} has a closed form solution, which we find
by taking the gradient with respect to each $\sigma_i^2$ and setting 
it to $0$: 
\begin{equation}
\label{OptimalSigmas}
\bar\sigma_i^2(x)  = \frac{1}{N_i} \|d_{i} - F_{i}(x)\|_2^2
\end{equation}
This quantity is precisely the population variance estimate. The modified problem~\eref{ProjectedClass}
is now given by 
\begin{equation}
\label{MAPmultipleMod}
\min_{x}\quad \tilde g(x) :=   
\sum_{i=1}^M \left( 
N_i\log(2\pi \bar\sigma^2_i(x)) 
+ 
N_i
\right)\;.
\end{equation}
The gradient of this objective is given by 
\begin{equation}
\label{GradMultVar}
\nabla_x \tilde g(x) = -\sum_{i=1}^m\frac{1}{\bar\sigma^2_i(x)}\nabla F_{i}(x)
(d_{i} - F_{i}(x))\;,
\end{equation}
while the Gauss-Newton (GN) Hessian approximation is given by 
\[
H(x) = \sum_{i=1}^m\frac{1}{\bar\sigma^2_i(x)}\nabla F_{i}(x) \nabla F_{i}(x)^T \;.
\]
Note that this is an approximation to $\nabla^2_x g$, and completely ignores 
the term $\nabla^2_{x,\theta} g(\bar x, \bar\theta(\bar x))\nabla_x\bar\theta(\bar x)$ in~\eref{secondProjDer}. 
The term can actually be explicitly calculated for this application, and turns out to be a dense negative definite 
correction to the Hessian approximation. If we ignore it, we 
recover the algorithm in~\cite{BellBurkeSchumitzky}, which can be seen by forming the GN subproblem:
\begin{equation}
\label{GNfunction}
\min_x  \sum_{i=1}^M \frac{1}{\bar\sigma^2_i(x_k)}
\|
d_{i} - F_{i}(x_k) - \nabla F_{i}(x_k)^Tx
\|_2^2\;.
\end{equation}
This expression matches~\cite[(12)]{BellBurkeSchumitzky} up to a constant.  
However, while in~\cite{BellBurkeSchumitzky} the function~\eref{GNfunction} came about as a cleverly 
constructed proxy objective for~\eref{MAPmultiple}, we can now view it as a natural GN approximation to the 
modified objective~\eref{MAPmultipleMod}. 

\paragraph{Example: Full Waveform Inversion}

Full waveform inversion (FWI) is an approach to obtain gridded subsurface velocity 
parameters from seismic data. Experiments are conducted by placing explosive
sources on the surface and recording the reflected waves with an array of
receivers on the surface. FWI is naturally cast as a nonlinear least squares
optimization problem \cite{Tarantola1984,Pratt1998}, and fits in the framework
described above. The data, $d_{i}$, in this case represents the Fourier
transform of the recorded time series for frequency $i$. The
corresponding modeling operator, $F_{i}(x) = PA_i(x)^{-1}Q_i$, inverts a discretized Helmholtz 
operator $A_i(x)$ for the $i^{\mathrm{th}}$ frequency and the gridded velocity field $x$ 
and samples the wavefield at the receiver locations. 
Here, $P$ denotes the sampling operator and each column of the matrix $Q_i$ is a gridded source function.

\begin{figure}
\centering
\begin{tabular}{ccc}
\includegraphics[scale=.3]{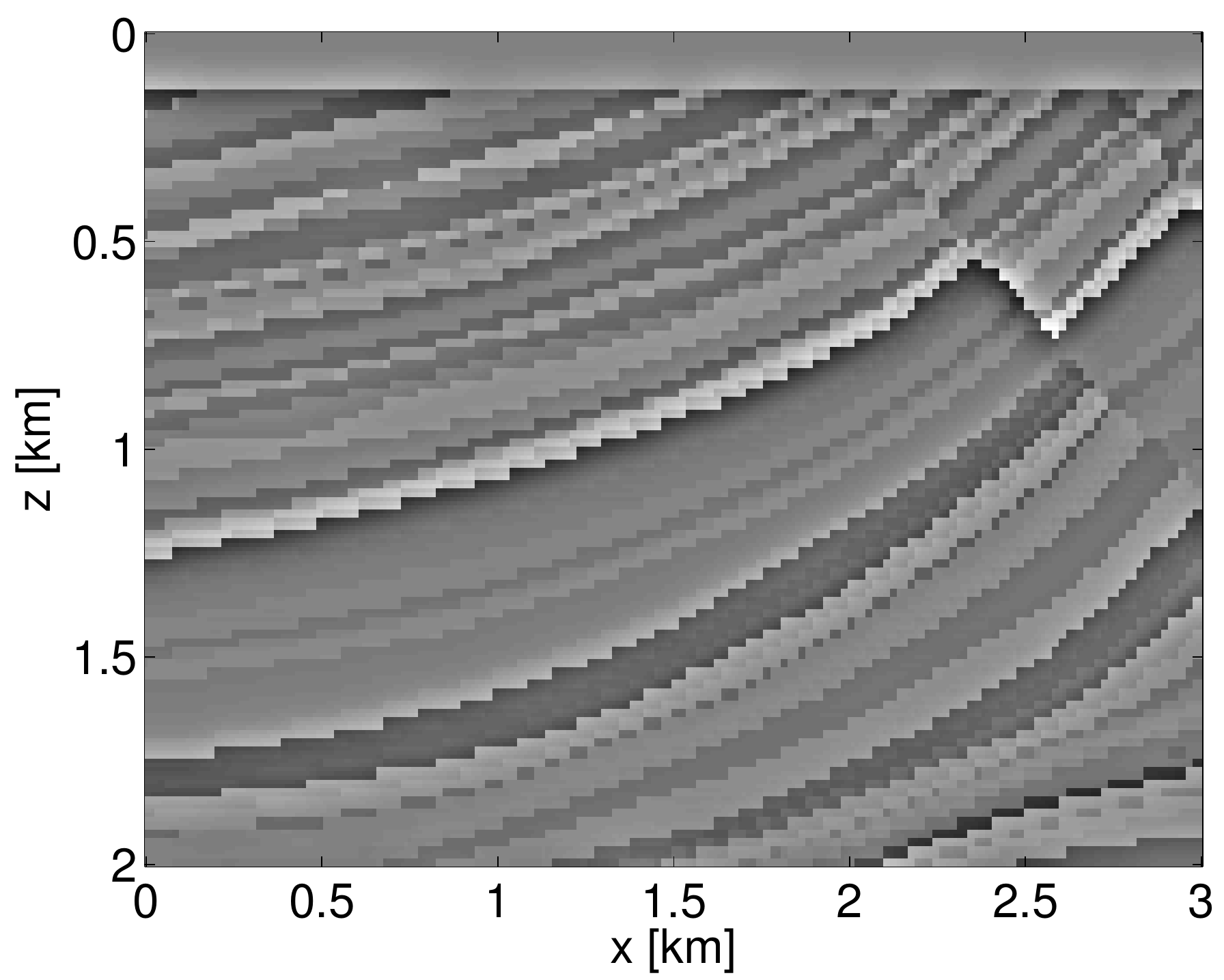}&
\includegraphics[scale=.3]{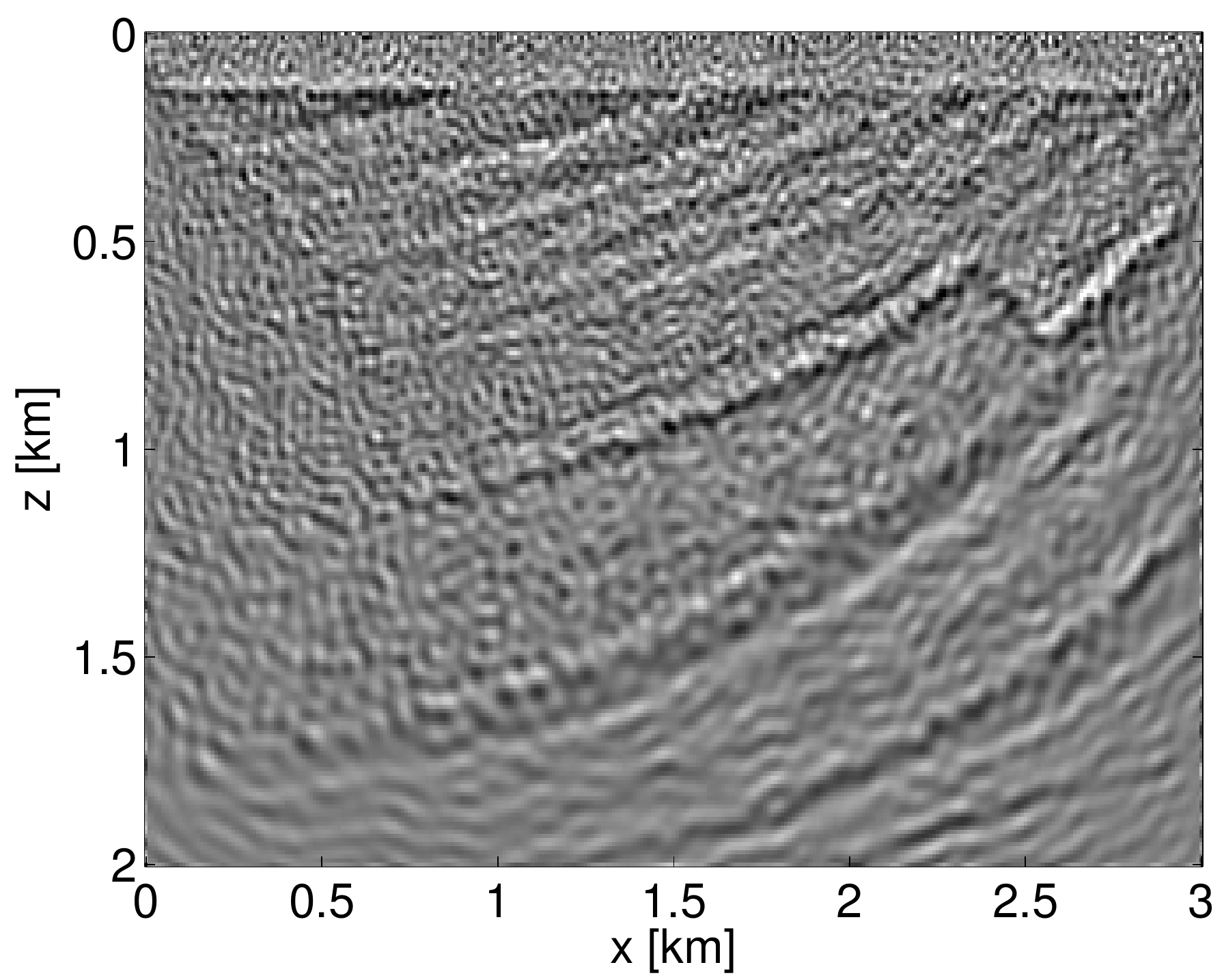}&
\includegraphics[scale=.3]{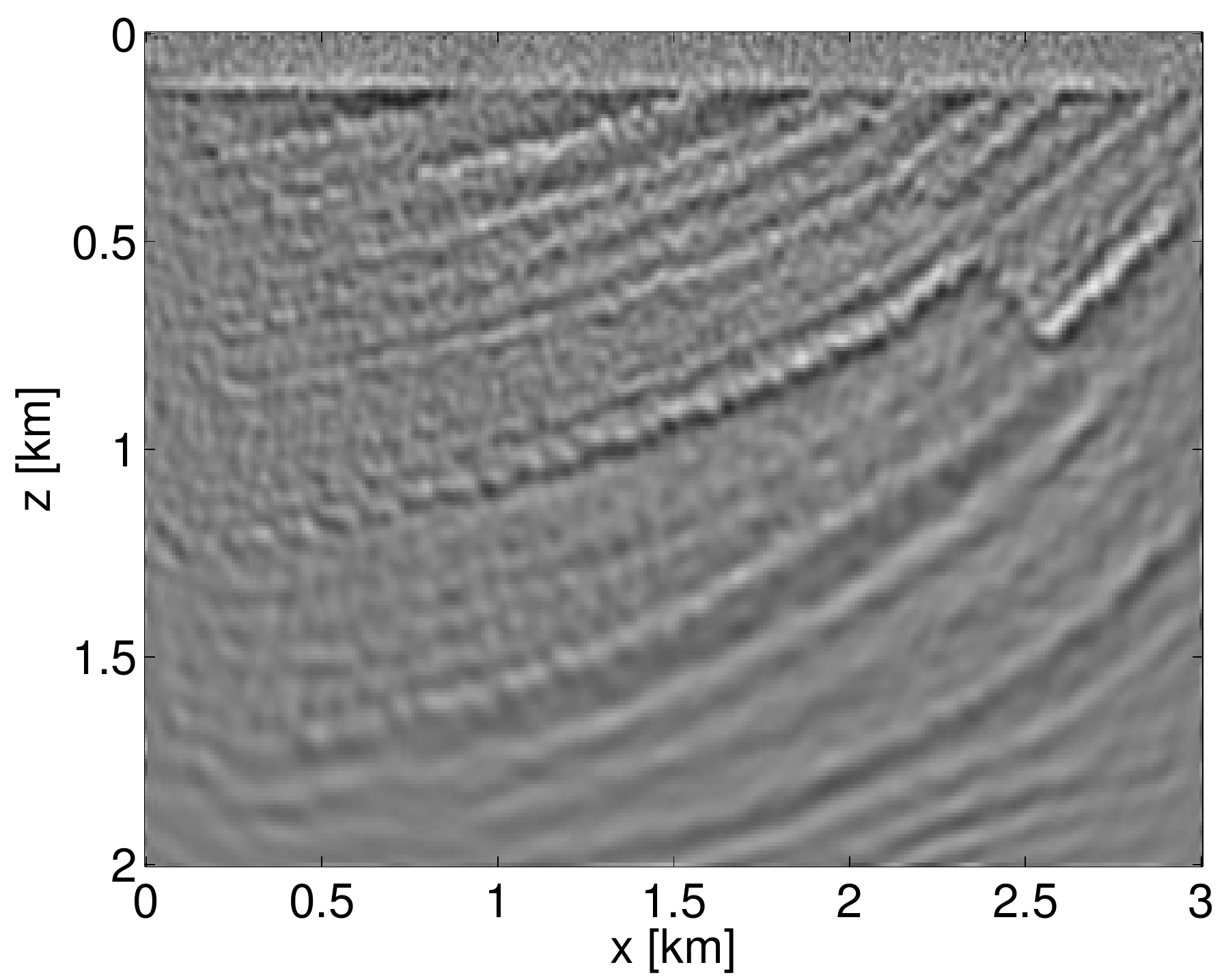}\\
{\small (a)}&{\small (b)}&{\small (c)}\\
\end{tabular}
\centering
\begin{tabular}{cc}
\includegraphics[scale=.3]{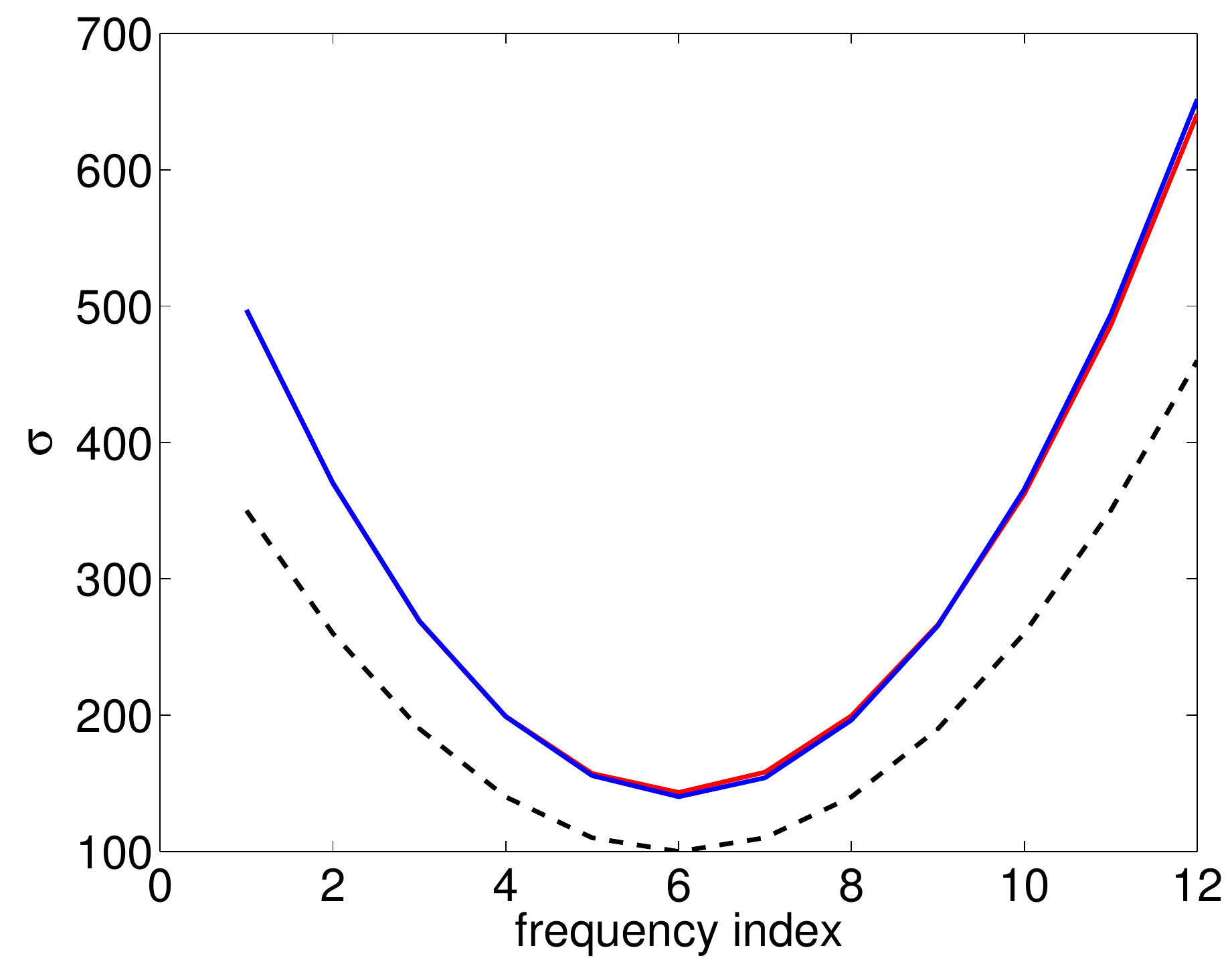}&
\includegraphics[scale=.3]{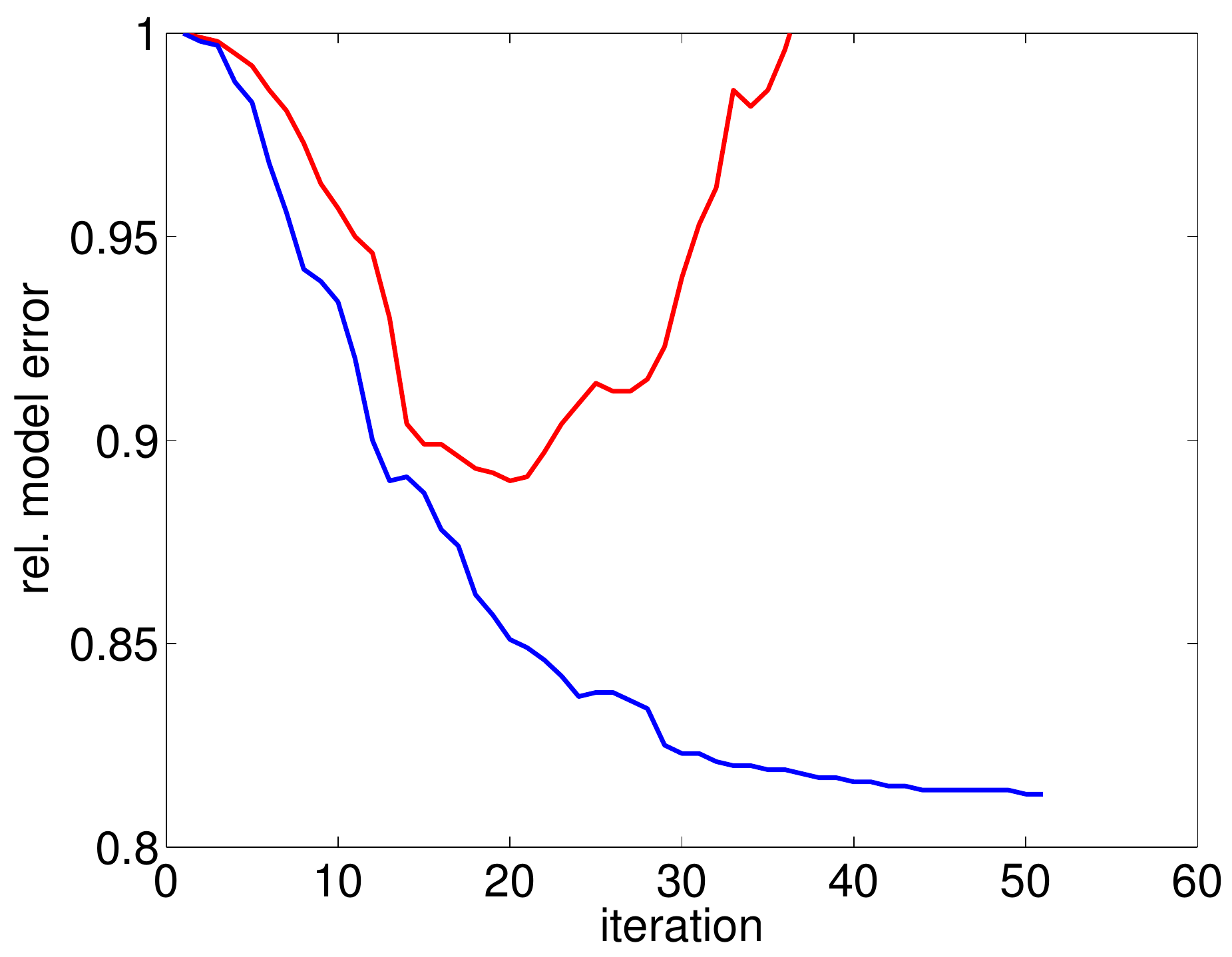}\\
{\small (d)}&{\small (e)}\\
\end{tabular}
\caption{Results for variance estimation. (a) True model, (b) result after 50
iterations for fixed constant variance $\sigma_i = 1$ and (c) result after 50 iterations with
variance estimation. The sample-variance for the latter two (red and blue respectively) results as well as
the true variance (dashed line) is shown in (d). Finally, (e) shows the relative model error
for each iteration for fixed (red) and estimated (blue) variance.}
\label{fig:varest}
\end{figure}

To illustrate the approach we use a subset of the well-known Marmousi benchmark
model, depicted in figure~\ref{fig:varest} (a). The model is discretized on a
201 $\times$ 301 grid with 10 m grid spacing. We generate data for 151 sources, 301 receivers
(i.e., $N_i = 151 \times 301$)---all equi-spaced and located at the surface---
and $M=12$ frequencies between 3 and 25 Hz. Typically, the data has a lower
signal to noise ratio for the low and high frequencies. To emulate this
situation we add Gaussian noise to the measurements with variance $\sigma_i \sim (i - 6)^2$. 
We use an L-BFGS method to  solve both the the modified optimization
problem~\eref{MAPmultipleMod} and the original problem for a fixed
$\sigma_i=1$ for all $i$. The results after 50 iterations are shown in
figure~\ref{fig:varest}(b,c). The
corresponding error between the reconstructed and true model is shown in
figure~\ref{fig:varest}(d). Finally, we show the estimated variance at the final model for both
reconstructions in figure~\ref{fig:varest}(e). The reconstruction obtained by solving the modified problem is
clearly better. Interestingly enough, the variance estimates for both models are almost
identical.

\subsection{Correlated Multivariate Observations}
\label{CorrelatedObservations}

The results from the previous case can be generalized to general variance estimation 
in a multivariate inverse problem setting with correlated errors. Consider the model~\eref{statisticalModel}, 
where now we take
\(\epsilon_i \sim N(0, \Sigma)\).
In this case, all of the $\epsilon_i$ are in of the same dimension. The ML objective corresponding to~\eref{MAPmultiple}
is given by 
\begin{equation}
\label{FullMAPmultiple}
\hspace{-1cm}
\min_{\Sigma, x} \quad
g(x, \Sigma) :=  
\left(M\log(2\pi \det(\Sigma)) 
+ \sum_{i=1}^M (d_i - F_i(x))^T \Sigma^{-1}(d_i - F_i(x))
\right)\;.
\end{equation}


The point here is that despite the generalization to full $\Sigma$, we still have a closed form solution
analogous to~\eref{OptimalSigmas}:
\begin{equation}
\label{OptimalSigmaFull}
\Sigma(x) = \argmin_\Sigma g(x,\Sigma) = \frac{1}{M}\sum_{i = 1}^M (d_i - F_i(x))(d_i - F_i(x))^T\;.
\end{equation}
This can be shown by a simple derivative computation: 
\[
\begin{array}{lll}
\frac{d}{d\Sigma^{-1}} g(x, \Sigma) &=& -M\Sigma + \frac{d}{d\Sigma^{-1}}\tr(\sum_{i=1}^m\Sigma^{-1}(d_i - F_i(x))(d_i - F_i(x))^T)\\
&=&-M\Sigma + \sum_{i=1}^m (d_i - F_i(x))(d_i - F_i(x))^T = 0\;.
\end{array}
\]
Therefore, the variable projection method applies immediately to~\eref{FullMAPmultiple}, at the cost of computing, 
at each iteration in $x$, the standard multivariate variance estimate~\eref{OptimalSigmaFull}. If this cost is high
(i.e. if each $\epsilon_i$ has high dimension), there are still a number of strategies that make the proposal feasible. 
For example,~\eref{OptimalSigmaFull} can be computed approximately using a random subset of the residuals. 

In addition to improving the primary parameters, incorporating nuisance parameter estimation 
can be helpful to post-processing analysis such as uncertainty quantification. 
For example, the estimate of $\Sigma$ at the final solution can be used to estimate the posteriori covariance matrix
in the model space $\left(\nabla F \Sigma^{-1} \nabla F^T\right)^{-1}$ ~\cite{Flath:2011}.

\subsection{Degrees of Freedom and Variance Estimation for Student's t Formulation}
\label{MetaParamSt}
Many applications require robust
formulations to obtain reasonable results with noisy data or 
in cases where a portion of the data is unexplained by the forward model (e.g. in the 
presence of {\it coherent artifacts}). 
A useful way to derive these formulations is to begin with the statistical
model~\eref{statisticalModel}
where the noise term $\epsilon_i$ is modelled using a particular parametric density, and then formulate the 
{\it maximum a posteriori} (MAP) likelihood problem.
The least-squares formulation corresponds to a Gaussian assumption 
on $\epsilon_i$ (see section \ref{VarEst}), while assuming a Laplacian distribution
leads to a one-norm penalty on the data-misfit.

As shown by~\cite[Theorem 2.1]{AravkinFriedlanderHerrmannVanleeuwen:2011}, 
in cases where unexplained artifacts may be large or constitute a significant portion of the data, 
it is better to use {\it heavy-tailed} densities. A prime example is the Student's t, whose density is given by 
\begin{equation}
\label{StudentDensity}
\mathbf{p}(y, \sigma^2, k) = \frac{\Gamma((k+1)/2)}{\Gamma(k/2)\sqrt{\pi k \sigma^2}}
\left( 
1 + \frac{y^2}{\sigma^2 k}
\right)
^{-(k+1)/2}\;.
\end{equation}

This density was first successfully used in~\cite{Lange1989} in the data fitting context. 
The degrees of freedom parameter $k$ was seen as a tuning parameter, smoothly transitioning
between heavy-tailed and near-Gaussian behaviour; $k$ and $\sigma$ were fit using Expectation Maximization (EM) and scoring methods.
This density was also successfully used in the Kalman smoothing
context~\cite{Fahr1998}, where it was suggested that the EM algorithm can be used to fit meta-parameters.  
Recent work using the Student's t distribution~\cite{AravkinFriedlanderHerrmannVanleeuwen:2011,RobustBAArxiv, aravkin12}
has side-stepped the problem, using fixed values for $\sigma$ and $k$.

In this section, we show that the general projection approach can be used to solve the 
joint inverse problem, treating scale and degrees of freedom as nuisance parameters. 
We propose a novel simple method, different from EM or scoring methods discussed in~\cite{Lange1989},
for estimating scale and degrees of freedom for any given set of residuals.
Given the model~\eref{statisticalModel}, the full MAP Student's t estimation problem is given by 
\begin{equation}
\label{FullStMap}
\hspace{-2cm}
\min_{x, k, \sigma^2} g(x, \sigma^2, k) 
:= 
-n\log\left(\frac{\Gamma\left(\frac{k+1}{2}\right)}{\Gamma\left(\frac{k}{2}\right)\sqrt{\pi k}}\right) + \frac{n}{2}\log(\sigma^2) + \frac{k+1}{2}\sum_{i=1}^n\log\left(1 + \frac{r_i^2}{\sigma^2 k}\right)\;,
\end{equation}
where $r_i = d_i - F_i(x)$. Following the philosophy presented in the paper,
we solve the problem by defining the modified objective
\begin{equation}
\label{StKS}
\tilde g(x) = g(x,\bar\sigma^2(x),\bar k(x))
\end{equation}
with 
\[
(\bar\sigma^2(x), \bar k(x)) = \argmin_{\sigma^2, k} g(x, \sigma^2, k)\;.
\]
The two-dimensional optimization problem in $(\sigma^2,k)$ required to evaluate $\tilde g(x)$ can be
solved using a customized routine or a black-box optimization code. An application is presented below.

\paragraph{Example: Traveltime tomography:}
We consider a cross-well traveltime tomography problem. In this case, sources
and receivers are placed in vertical wells and the data consists of 
picked traveltime of first arrivals. Since the data are typically very noisy,
a portion of the traveltimes may be picked erroneously, motivating the use of
robust penalties for the inversion. The traveltimes are computed by a geometric
optics approach, where wave propagation is modeled via rays. The traveltime
between a given source and receiver is simply the integral of the reciprocal
velocity along the corresponding ray-path. By assuming small perturbations of a
known \emph{background} velocity, we arrive at a linear modeling operator with a
fixed ray geometry. The data are the traveltime perturbations,
while the primary parameter of interest is the velocity perturbation, both taken with respect to 
a known background model.

\begin{figure}
\centering
\begin{tabular}{cc}
\includegraphics[scale=.3]{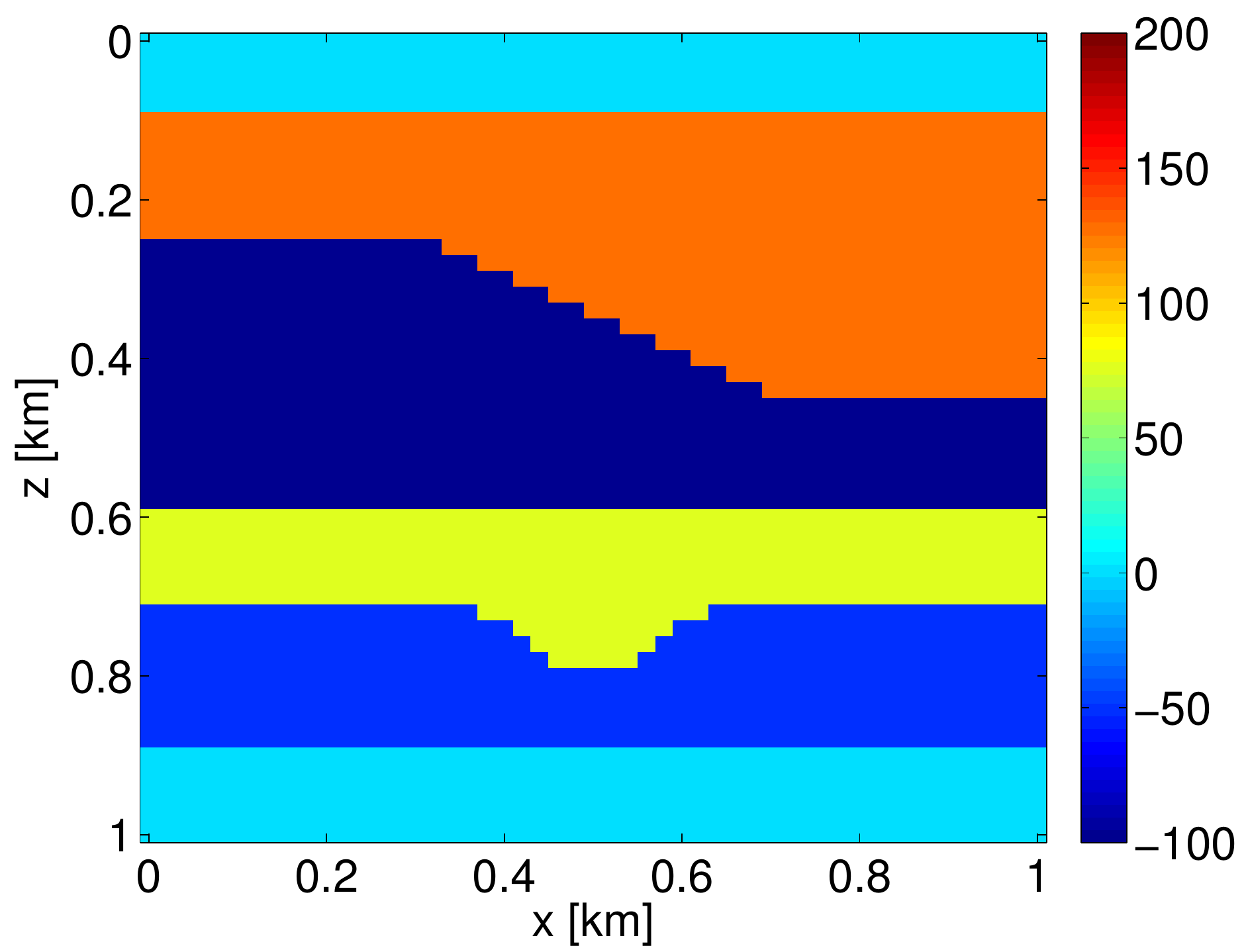}&
\includegraphics[scale=.3]{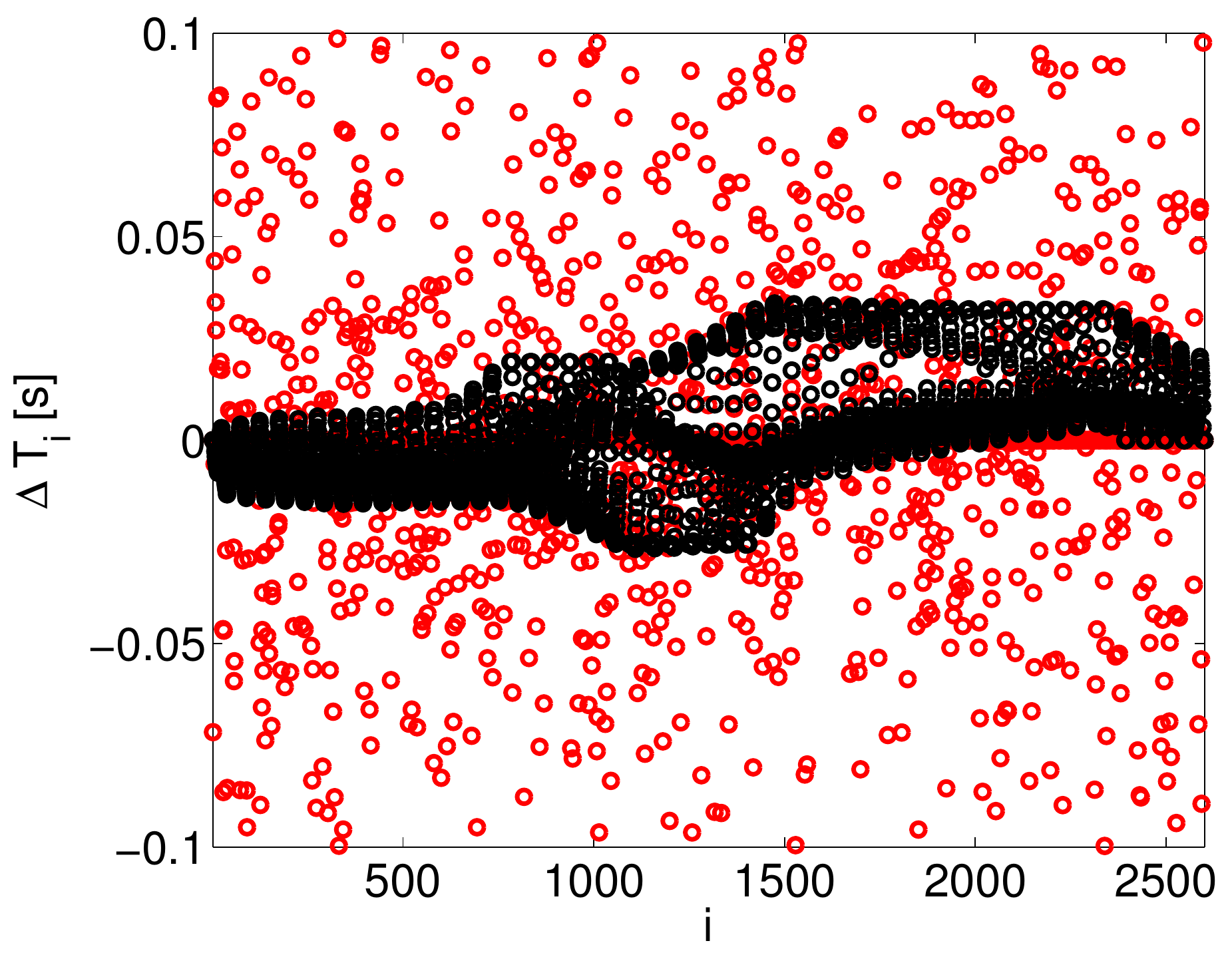}\\
{\small (a)}&{\small (b)}\\
\end{tabular}
\caption{(a) velocity perturbation in m/s used to generate the observed data, (b)
shows the corresponding traveltime perturbations in black and the outliers in
red.}
\label{fig:dfest1}
\end{figure}

In this example, we consider a constant background velocity, so that the ray
paths are straight lines. The modeling operator is therefore essentially a Radon
transform,  which is often used in medical X-ray imaging applications. The true
velocity perturbation is discretized on a 51 $\times$ 51 grid and  is shown in figure~\ref{fig:dfest1}(a).
The corresponding data for 51 sources and receivers and the added outliers are shown in figure~\ref{fig:dfest1}(b).
We regularize the inversion by inverting for the primary parameters on a courser grid of
size 26 $\times$ 26. We then interpolate back to the fine grid using 2D cubic
interpolation.
The modified optimization problem is now given by
\begin{equation}
\label{TT}
\min_x \rho_{\theta}(\Delta T - Ax),
\end{equation}
where $\Delta T \in \mathbb{R}^{2601}$ are the measured traveltime
perturbations, $x\in\mathbb{R}^{676}$ is the velocity
perturbation, $A$ is the modeling operator which combines the Radon transform and
interpolation, and $\theta = (\bar\sigma^2, \bar k)$
is obtained by solving~\eref{StKS} using a Nelder-Mead method~\cite{nelder}.

Note that we may treat $\theta$ as fixed when
designing an algorithm to solve~\eref{TT}, as long as the parameters are
re-estimated at every evaluation of $\rho_{\theta}(r)$ and its derivatives. 
To solve~\eref{TT} we use a modified Gauss-Newton algorithm which
calculates the updates by solving
\begin{equation}
\left(A^TH_{\theta}(r_k)A\right)\Delta x_k = A^T\nabla\rho_{\theta}(r_k),
\end{equation}
where $r_k = \Delta T - Ax_k$ and $H_{\theta}$ is a positive approximation of
the Hessian $\nabla^2\rho_{\theta}$. 
We solve the subproblems using CG. Note that when $\rho = \|\cdot\|^2$, 
the algorithm converges in one GN iteration, which is computed by 
applying CG to the normal equations. 

\begin{figure}
\centering
\begin{tabular}{ccc}
\includegraphics[scale=.3]{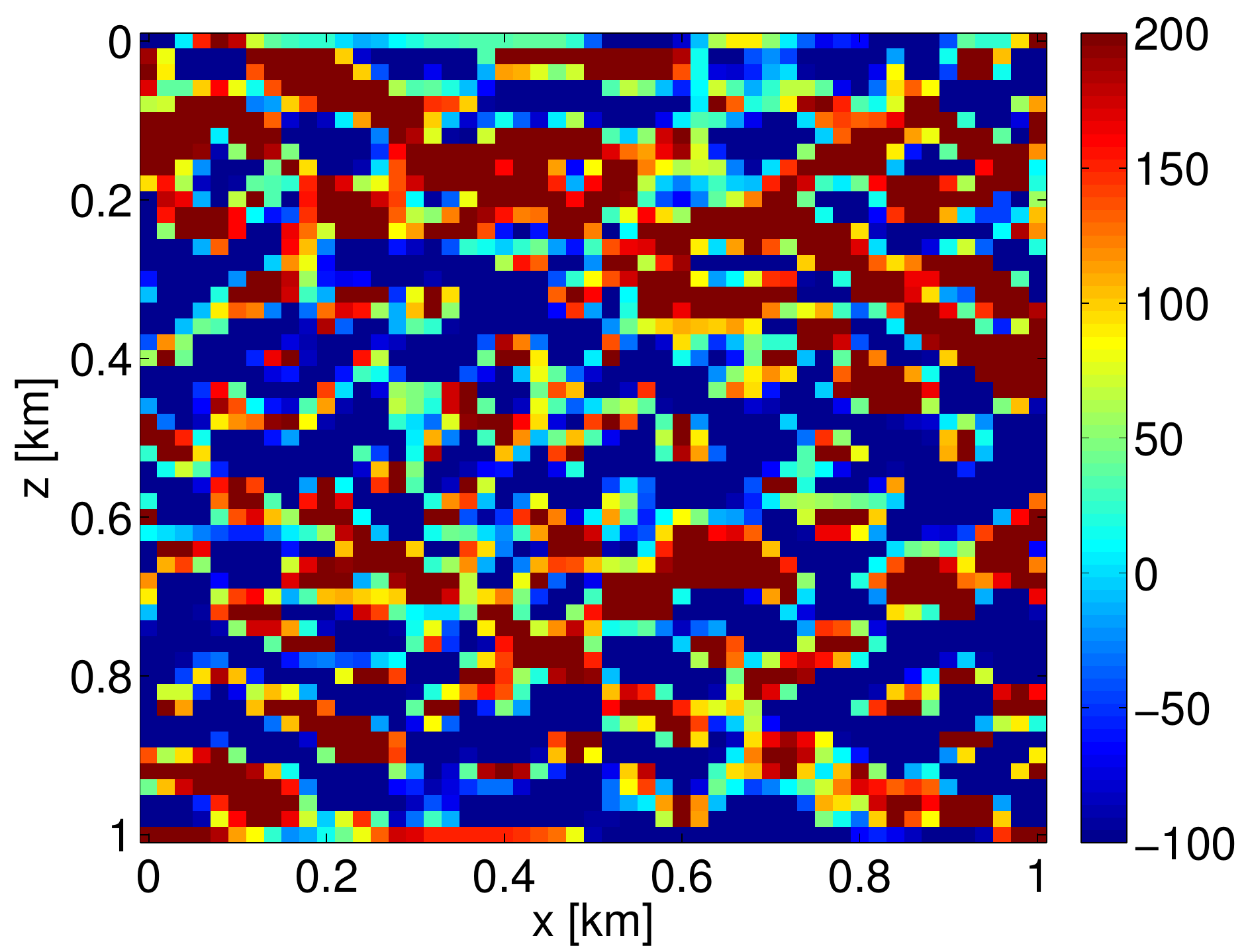}&
\includegraphics[scale=.3]{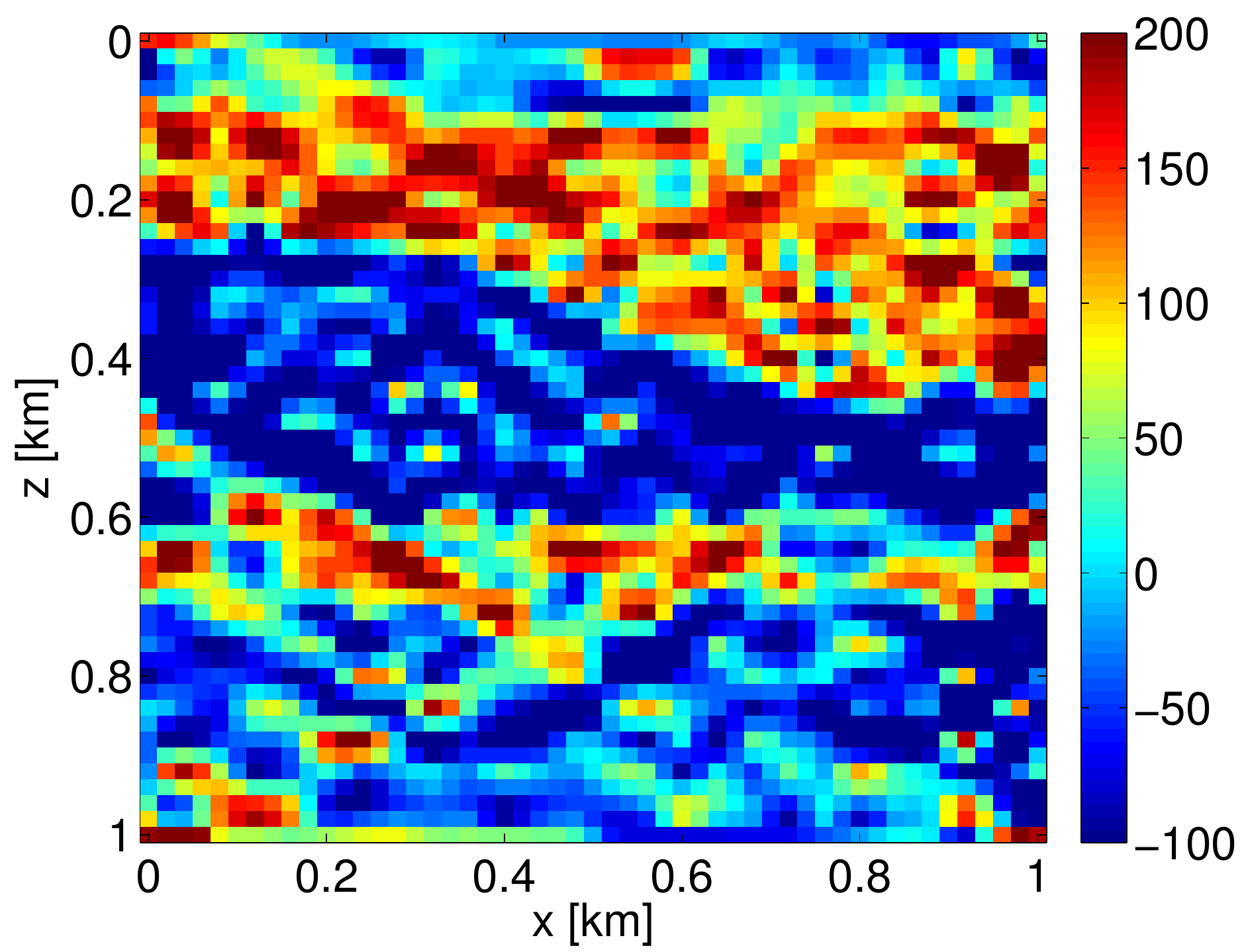}&
\includegraphics[scale=.3]{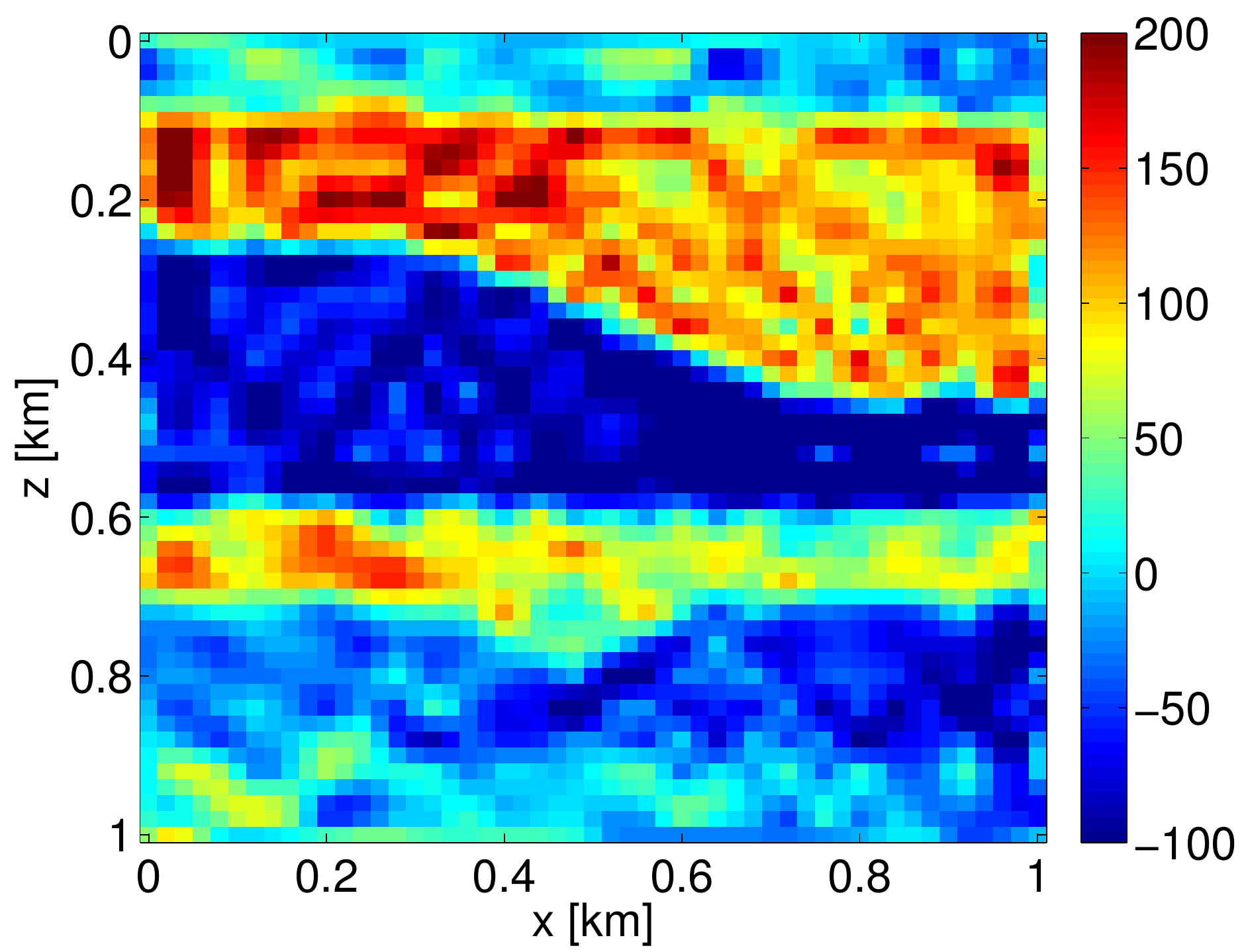}\\
{\small (a)}&{\small (b)}&{\small (c)}\\
\end{tabular}
\caption{Results for traveltime tomography. (a) least-squares reconstruction, (b)
Student's t reconstruction with fixed $\theta$ estimated at the initial residual
and (c) Student's t reconstruction where $\theta$ is re-estimated at every
iteration.}
\label{fig:dfest2}
\end{figure}

\begin{figure}
\centering
\begin{tabular}{ccc}
\includegraphics[scale=.3]{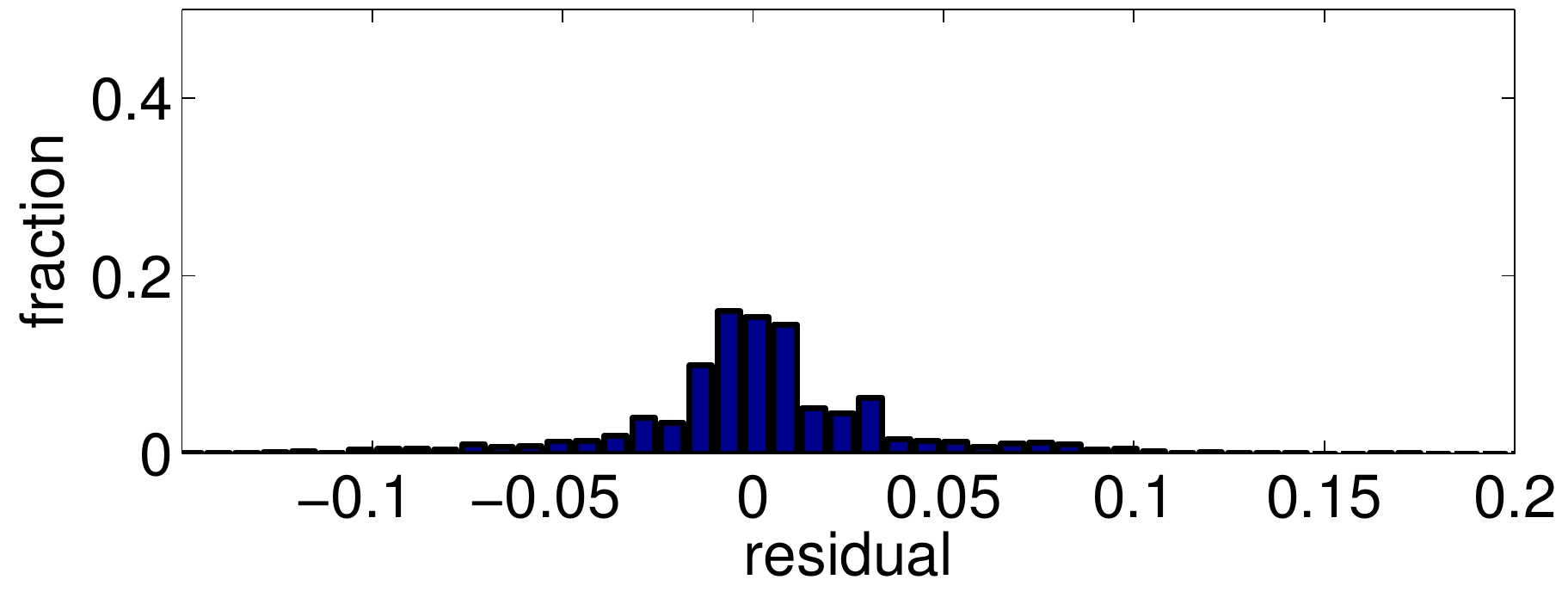}&
\includegraphics[scale=.3]{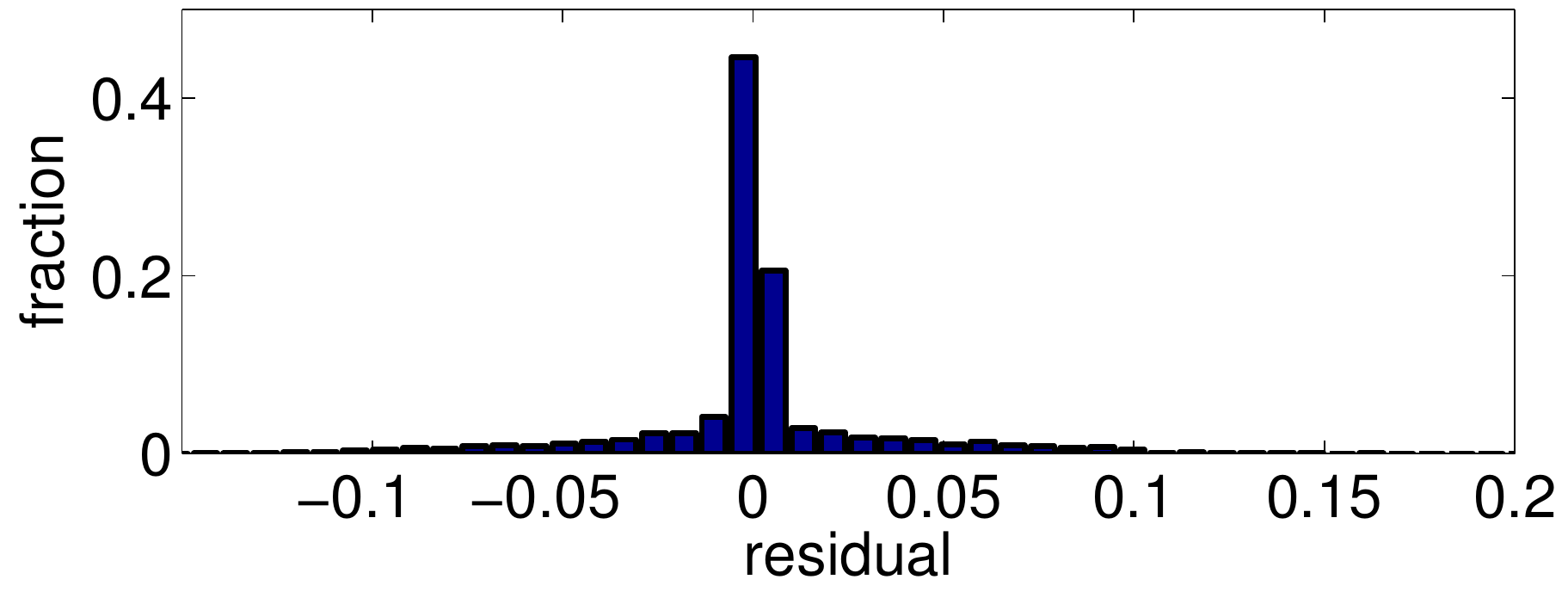}&
\includegraphics[scale=.3]{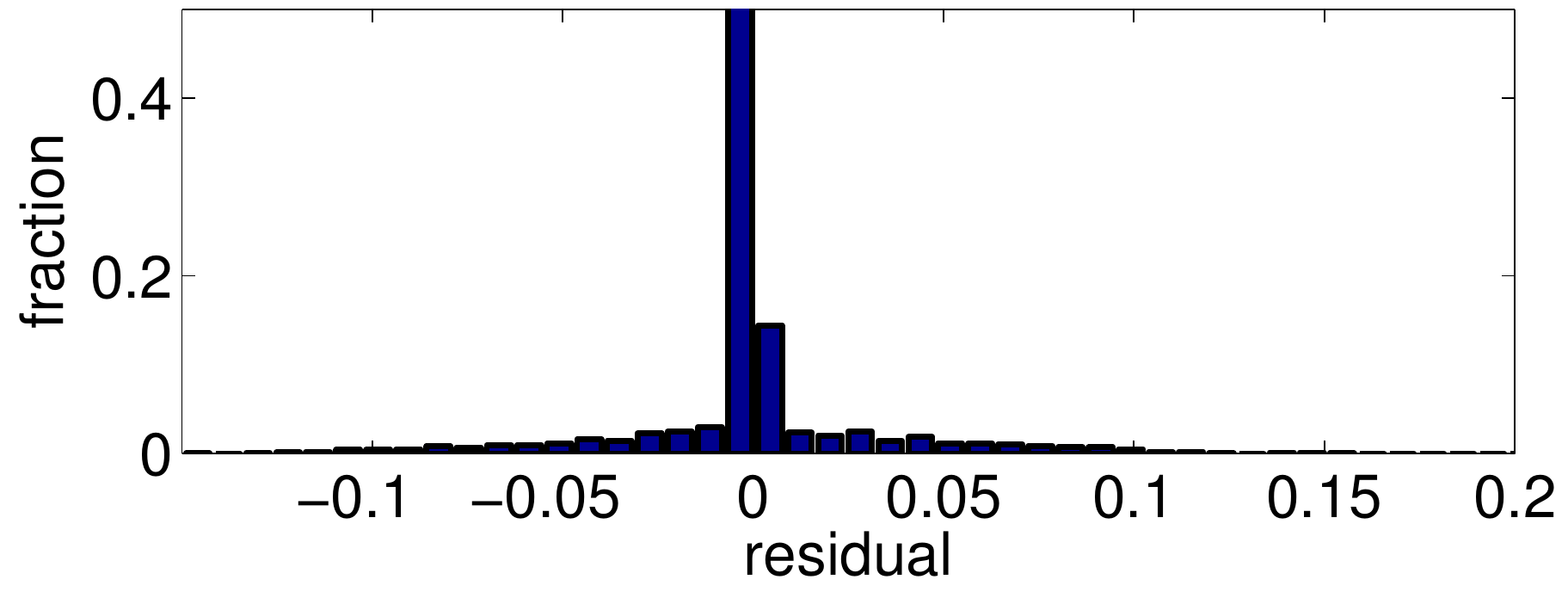}\\
\includegraphics[scale=.3]{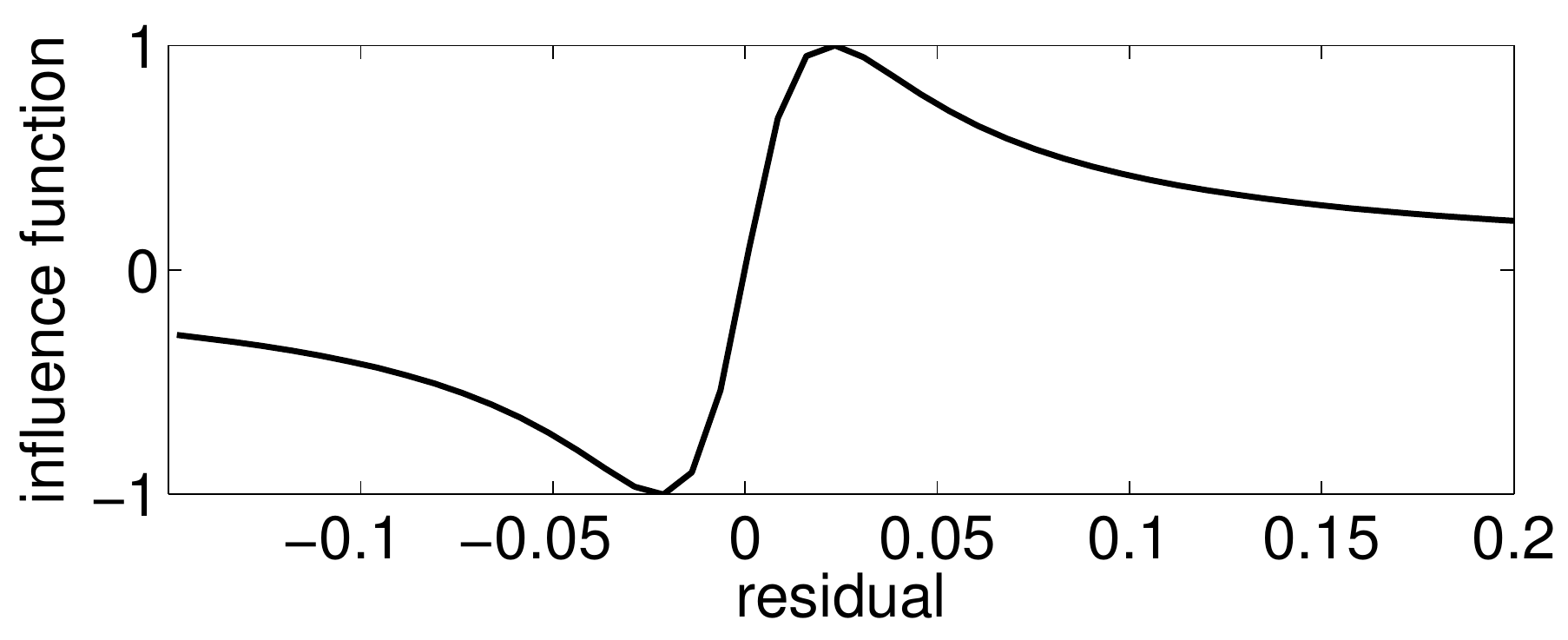}&
\includegraphics[scale=.3]{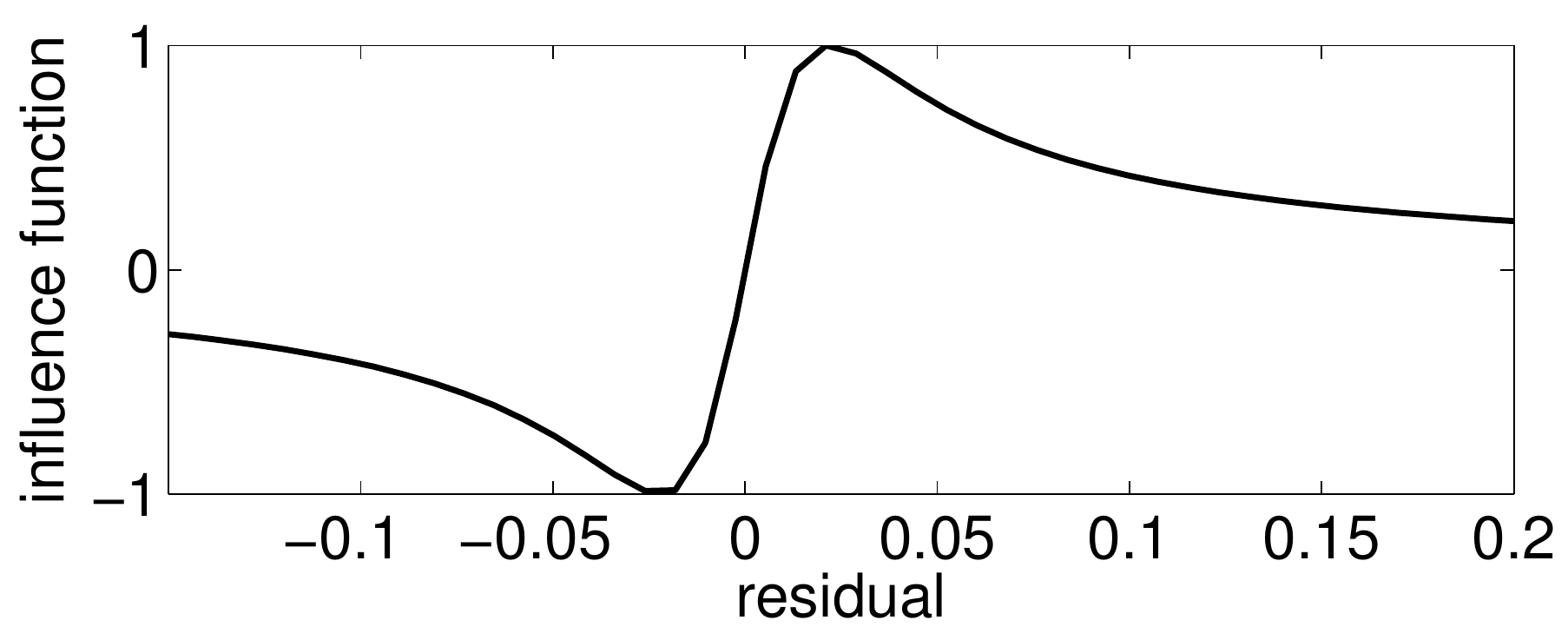}&
\includegraphics[scale=.3]{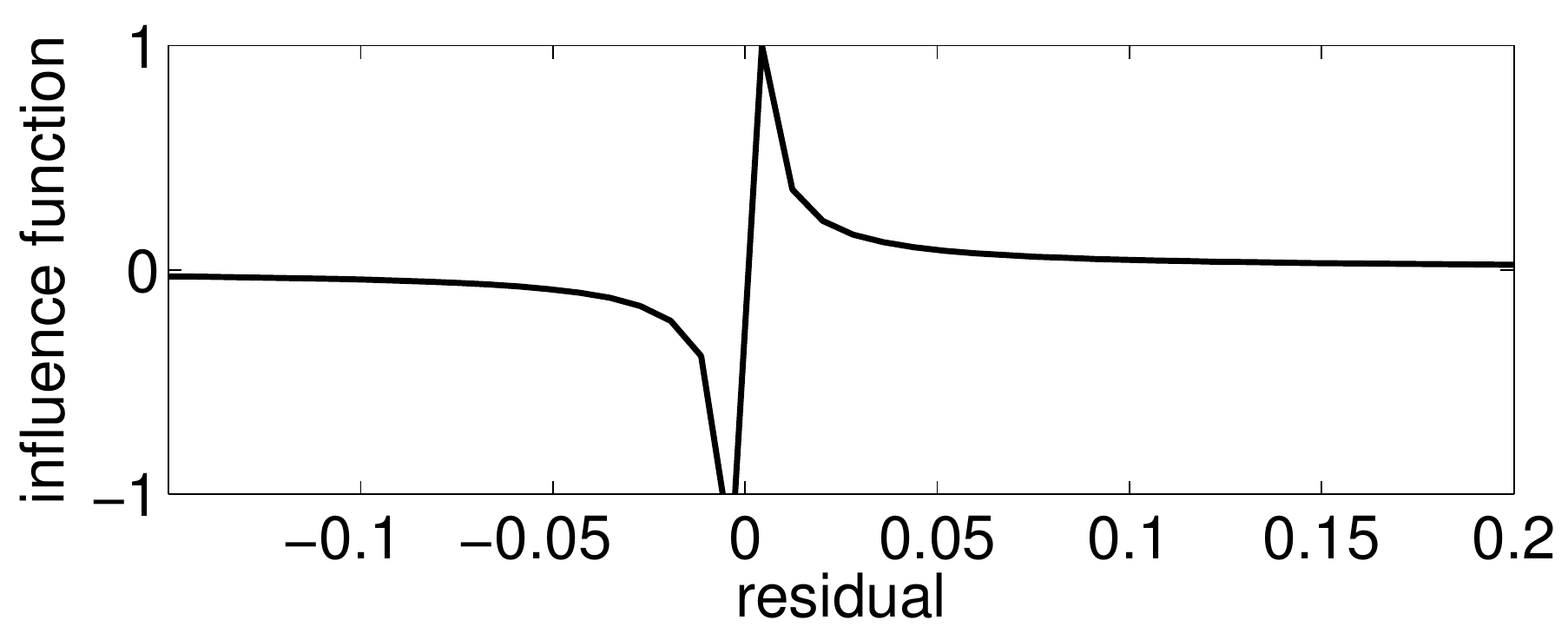}\\
{\small (a)}&{\small (b)}&{\small (c)}\\
\end{tabular}
\caption{Histograms of the residuals and corresponding influence functions
$\rho'_{\theta}$ . (a) initial residual, (b) final residual corresponding to
figure~\ref{fig:dfest2} (b) and (c) final residual corresponding to
figure~\ref{fig:dfest2} (c). In the latter case the parameters $\theta$ are
re-estimated at each iteration, allowing the inversion to home in on the good
data and ignore the outliers.}
\label{fig:dfest3}
\end{figure}

We compare the following three approaches: \emph{i)} least-squares,
shown in figure~\ref{fig:dfest2}(a), \emph{ii)} Student's t  with a fixed $\theta$ wich
is estimated once at the initial residual, shown in figure~\ref{fig:dfest2}(b),  and \emph{iii)} Student's t
where we estimate $\theta$ at each iteration, shown in
figure~\ref{fig:dfest2}(c). In order
to understand the difference between the latter two cases, we show histograms
of the initial and final residuals as well as the influence function for the
corresponding $\theta$ in figures~~\ref{fig:dfest3}(a-c). Clearly, re-fitting the $\theta$ at each iteration
allows the inversion to home in on the good data while ignoring the outliers.

\section{Automatic calibration}
\label{RobustSourceParam}
In this section we consider the case where the forward model includes a \emph{calibration factor} 
$\alpha$:
\begin{equation}
d = F(x,\alpha) + \epsilon.
\end{equation}

In case of the non-linear data-fitting problem
described earlier, the modified objective is given by 
\begin{equation}
\min_{x} 
\tilde g(x)=\rho(d- F(x,\bar\alpha)),
\end{equation}
where
\begin{equation}
\bar\alpha(x) = \argmin_\alpha\rho(d- F(x,\alpha)).
\end{equation}

The motivating example that led us to consider this class of problems is presented below. 

%
%

\paragraph{Example: FWI with source estimation:}
Seismic data can be interpreted as the Green's function of the subsurface,
parametrized by $x$, convolved with an unknown (bandlimited) source signature. In the frequency-domain
we can model the uknown source signature by multiplication with a complex scalar for each frequency.
The problem of interest is now formulated as
\begin{equation}
\min_{x,\alpha} \left\{ 
g(x,\alpha):=\sum_{i=1}^{M}\rho(d_i - \alpha_iF_i(x))
\right\}\;,
\end{equation}
where the index $i$ runs over frequency.
Just as in the variance parameter case, the parameters $\alpha_i$ are linked only through 
the parameters $x$, and for a given $x$ the problem decouples completely, giving
\begin{equation}
\label{AlphaOfx}
\bar\alpha_{i}(x) = \argmin_{\alpha_i}\rho(d_{i} - \alpha_iF_{i}(x))\;.
\end{equation}
We consider the least-squares and Student's t penalty, and use a
scalar Netwon-type method to solve~\eref{AlphaOfx}:
\begin{equation}
\alpha_i^{\nu+1} 
= 
\alpha_i^\nu 
- 
\frac
{\langle \nabla\rho(r_i^{\nu}(x)),F_i(x) \rangle}
{\langle F_i(x), HF_i(x)\rangle}\;,
\end{equation}
where $r_i^\nu(x) = d_i - \alpha_i^\nu F_i({x})$, $\nabla\rho$ is the gradient
of the penalty function and $H$ is (a positive definite approximation of) the Hessian
$\nabla^2\rho$. In particular, we have:

\begin{itemize}

\item{least-squares:} 
$\rho(r) = \frac{1}{2}\sum_i r_i^2$, 
$\nabla\rho_i = r_i$
and $H_{ii} = 1$.


\item{Student's t:} 
$\rho({r}) = \frac{1}{2} \sum_{i} \log(k + r_{i}^2)$,
$\nabla\rho_i = r_i/(k+r_i^2)$ 
and $H_{ii} = 1/(k+r_i^2)$.

\end{itemize}
For more details on the Student's t approach we refer to
~\cite{AravkinFriedlanderHerrmannVanleeuwen:2011}.

\begin{figure}
\centering
\begin{tabular}{cc}
\includegraphics[scale=.5]{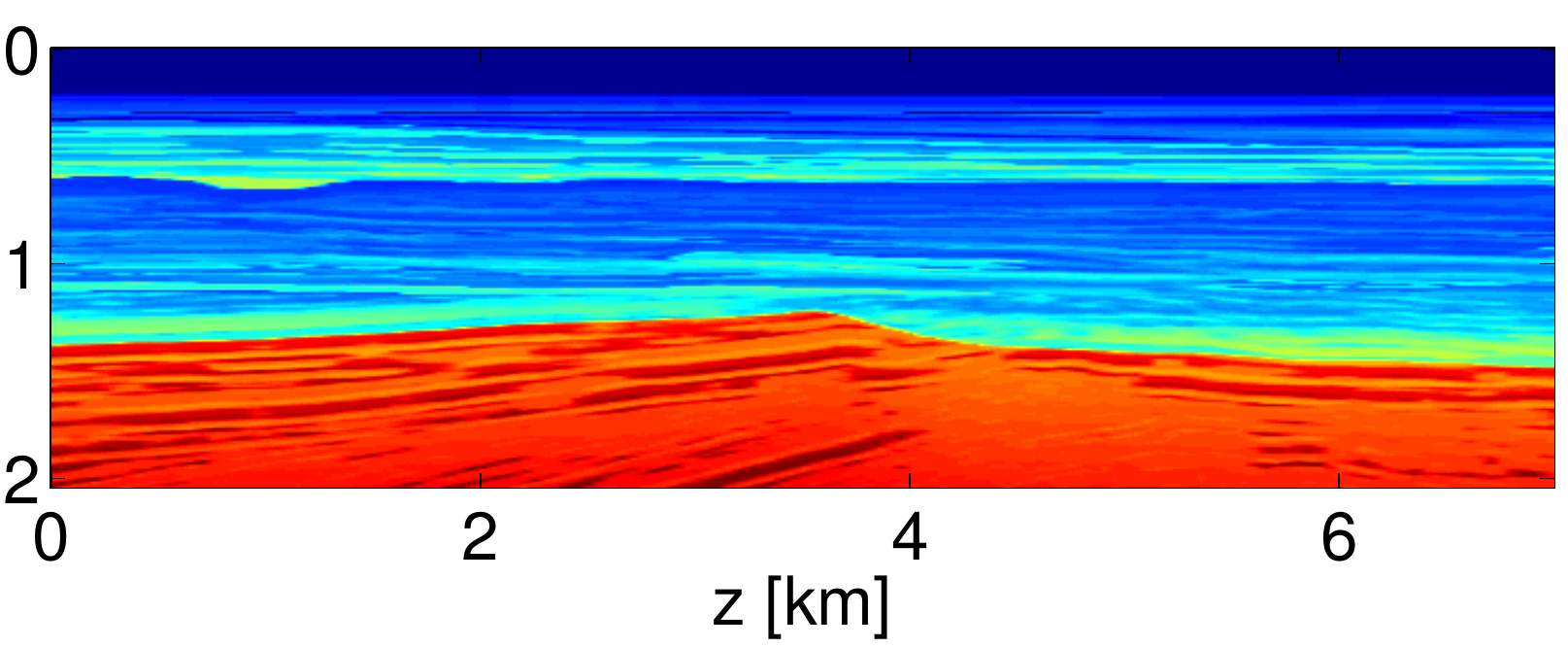}&
\includegraphics[scale=.5]{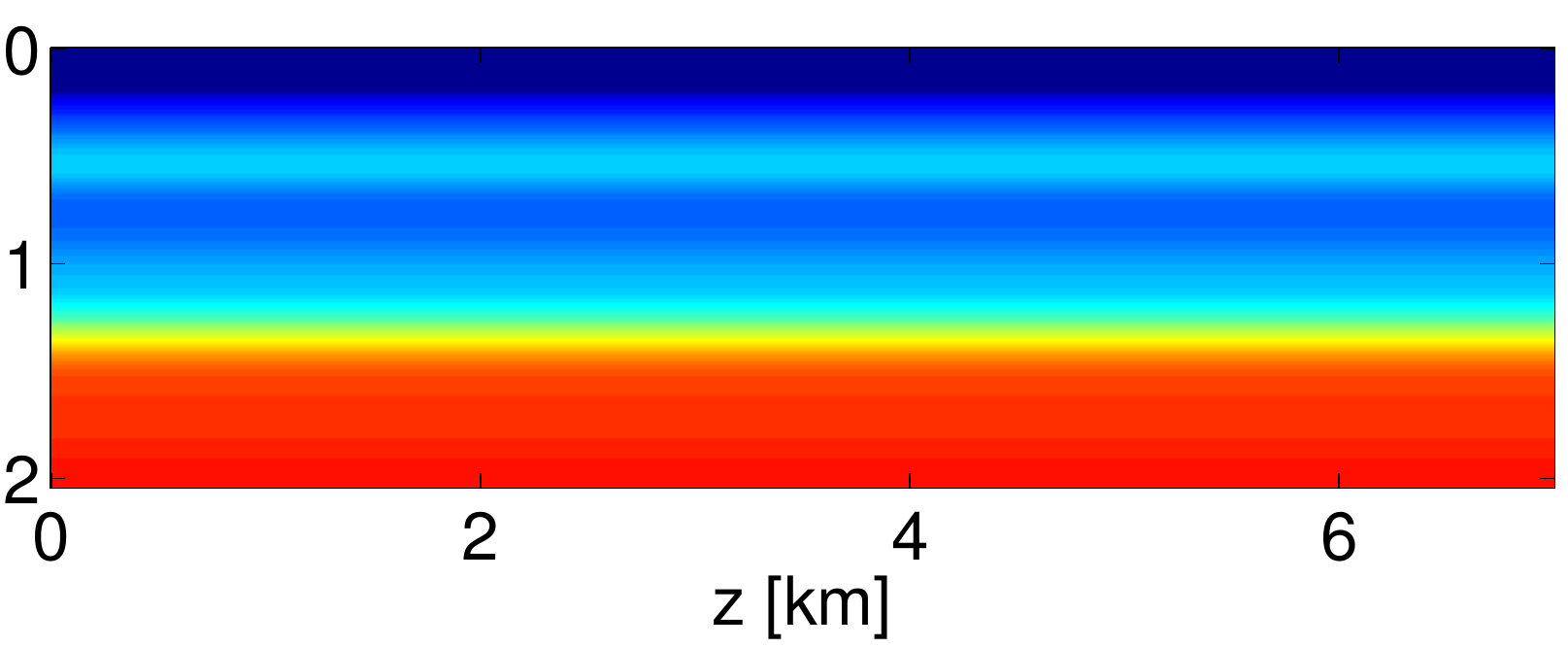}\\
{\small (a)}&{\small (b)}\\
\includegraphics[scale=.5]{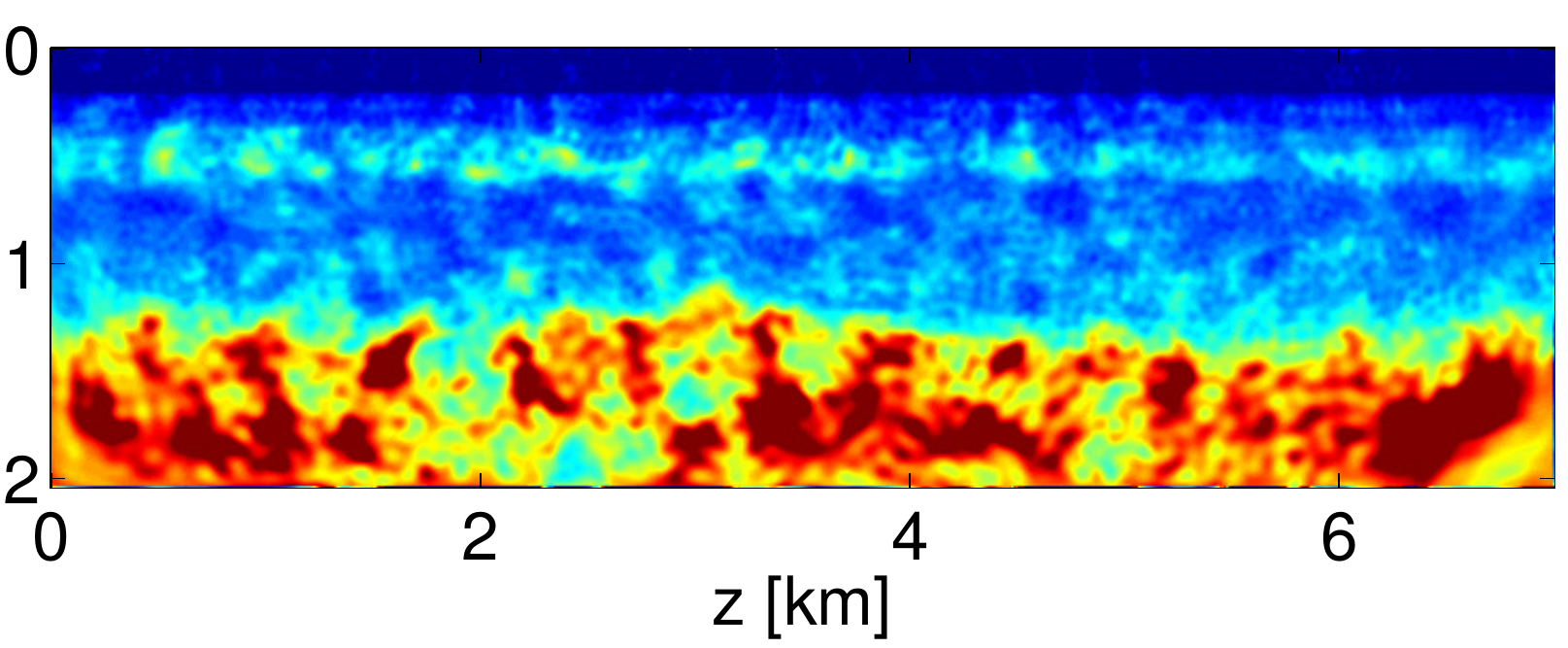}&
\includegraphics[scale=.5]{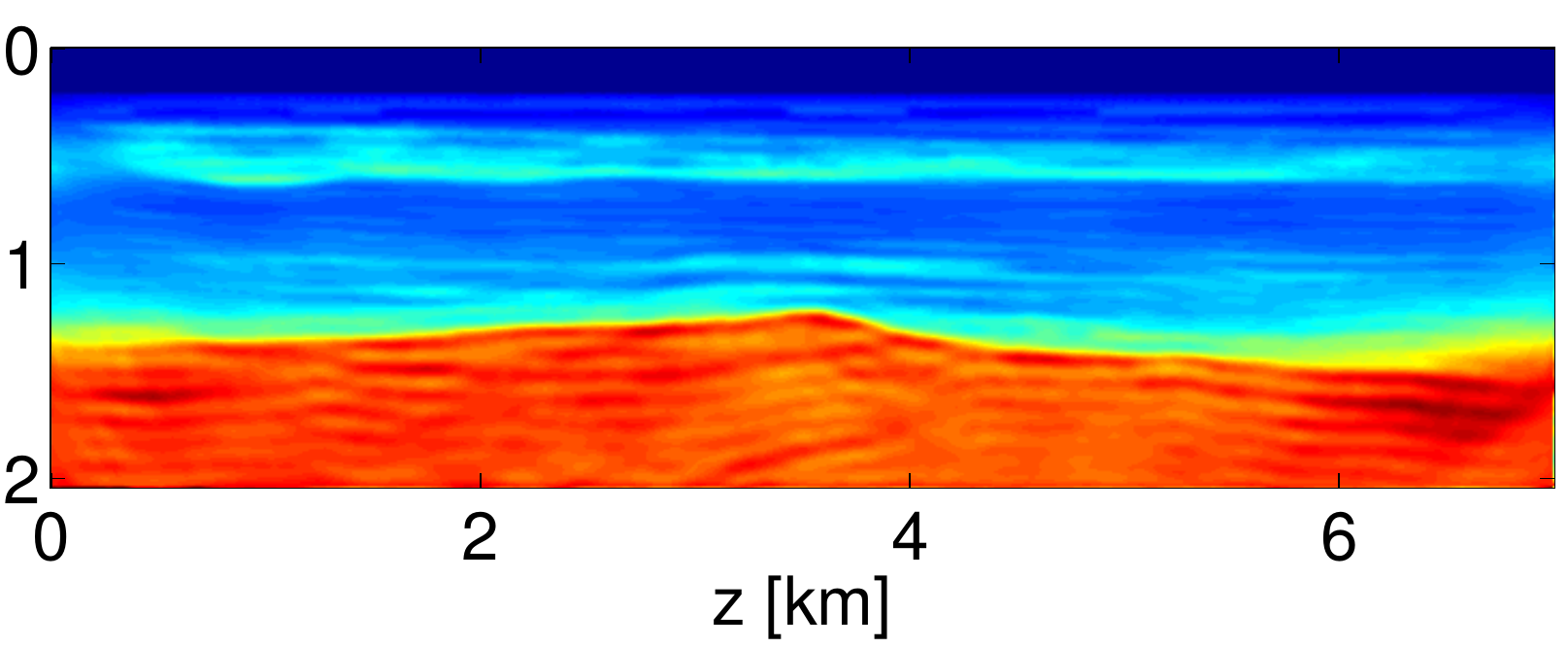}\\
{\small (c)}&{\small (d)}\\
\end{tabular}
\caption{(a) True velocity model. (b) Initial velocity model. (c)  Least-squares reconstruction from noisy data. (d) Student's t reconstruction from noisy data.} 
\label{fig:srcest}
\end{figure}

We generate seismic data for the velocity model depicted in figure \ref{fig:srcest} (a) with a time-domain
finite difference code. The data consists of 141 sources and 281 receivers and has a recording time of 4 seconds. 
10 percent of the data is corrupted with large outliers.

We invert the data in several stages, moving from low to high frequencies.   
Each stage uses only a few frequencies and the output is used as initial guess for the subsequent stage. 
This is a well-known strategy in FWI to avoid local minima~\cite{bunks95}.
We use an L-BFGS method to solve the resulting optimization problems, starting from the initial 
model shown in figure \ref{fig:srcest} (b). The
results are shown in figure \ref{fig:srcest} (c,d). The Student's t approach recovers the most
important features of the model whereas the least-squares approach leads to a very noisy model.

\section{Optimal experimental design}
\label{OED}
In Optimal Experimental Design one is concerned with finding \emph{optimal}  design 
parameters $\theta$ for which a set of test models $\{x_i\}$ can be recovered from the
corresponding simulated data $d_i(\theta) = F(x_i,\theta) + \epsilon$. This can be formulated as an 
optimization problem~(cf. \cite{horesh11} and references cited therein)
\begin{equation}
\label{OptDesign}
\min_{\theta} \sum_i Q(\bar{x}_i(\theta),x_i) + C(\theta),
\quad \bar{x}_i(\theta) = \argmin_{x} ||F(x,\theta) - d_i(\theta)||_2.
\end{equation}
Here, $Q(\bar{x}_i,x_i)$ measures the quality of the reconstruction (lower is better) and
$C(\theta)$ measures the \emph{cost} of a given experimental parameter setting. 
Note that~\eref{OptDesign} is actually the {\it reduced problem} for the joint optimization problem
\begin{equation}
\min_{x,\theta} \sum_i ||F(x,\theta) - d_i(\theta)||_2^2 + \sum_i Q(x,x_i) + C(\theta),
\end{equation}
where $x$ has been projected out. 

The downside to this approach is that projecting out $x$ is expensive, since {\it every iteration} 
of any algorithm to find $\bar x_i(\theta)$ requires repeated evaluations of $F(x, \theta)$, for {\it each model}
in the class $\{x_i\}$. 
Rather than projecting out $x$, we can project out the design parameter $\theta$ to arrive at 
a different reduced objective 
\begin{equation}
\min_{x} \tilde g(x) := \sum_i ||F(x,\bar\theta(x)) - d_i(\bar\theta(x))||_2^2 + \sum_i Q(x,x_i) + C(\bar\theta(x)),
\end{equation}
where $\bar\theta(x) = \argmin_{\theta} \sum_i ||F(x,\theta) - d_i(\theta)||_2^2 + C(\theta)$.
In many cases, $\bar \theta(x)$ can be computed cheaply without re-evaluating the whole forward model.
The forward modeling need only be done when $x$ is updated, as in the other applications that we presented. 
For example, $\theta$ may represent a vector of source weights for waveform inversion in which case
$F(x,\theta) = PA(x)^{-1}Q\mathsf{diag}(\theta)$ (see section \ref{VarEst}). 
Since the data are linear in the source, we need only invert the Helmholtz system once (for a given $x$), 
and therefore $\bar{\theta}$ can be computed relatively cheaply. 
A reasonable penalty on $\theta$ might be the one-norm, 
in which case we are looking for a setup with as few sources as possible. 
Alternatively, we can impose a two-norm penalty on $\theta$ to find a setup where the 
sources require the least amount of energy.

\section{Discussion and Conclusions}

Many inverse problems involve nuisance parameters that are not of primary interest
but can have significant influence on the estimation of primary parameters. Common examples
include variance, degree of freedom, and calibration parameters. 
These issues arise in a great variety of applications, including pharmacokinetic modeling~\cite{BellBurkeSchumitzky}, 
seismic inverse problems~\cite{Tarantola1984}, dynamic systems~\cite{fahrmeir2010multivariate}, 
uncertainty quantification~\cite{Flath:2011} and optimal experimental design~\cite{horesh11}. 

In this paper, we proposed a straightforward approach to fitting these nuisance parameters
on the fly, while solving the overall inverse problem. Specifically, we formulated the problem
as a joint optimization over primary parameters $x$ and nuisance parameters $\theta$~\eref{InverseProblemClass}, 
and showed that for a large class of problems, one can simply project out the $\theta$ parameters by 
solving~\eref{Projection}.  In this least squares case, this idea has been carefully studied under 
the name {\it variable projection}~\cite{Osborne2007,GolubPereyra2003}. 
As we showed, these ideas extend nicely to the entire class~\eref{InverseProblemClass}. 
In particular, Theorem~\ref{ProjectionTheorem} and Corollary~\ref{ConstrainedCorollary} characterize the general approach
and are the basis for algorithm design of first and second order methods. 

An immediate consequence of the work is the ability to modify first and second order algorithms 
that exploit particular application structure to also fit nuisance parameters. We 
 demonstrated this in practice using several (large scale) inverse problems:\\

\begin{tabular}{|c|c|c|}\hline
{\it Application} & {\it Complicating parameters} & {\it Algorithm}\\ \hline
full waveform inversion & variances in multiple datasets & L-BFGS\\\hline
travel time tomography & student's t parameters & Gauss-Newton with CG \\ \hline
automatic calibration & unknown source amplitudes & L-BFGS\\ \hline
\end{tabular}\\


In the case of variances in multiple datasets, the proposed approach matches the algorithm 
proposed in~\cite{BellBurkeSchumitzky}, and therefore the development we presented 
provides an alternative (and significantly simpler) derivation. We have also shown 
that the approach can be easily extended to estimate covariances between error sources
in the case where we have multivariate observations in Section~\ref{CorrelatedObservations}. 

In the case of student's t parameters, it is interesting to note that 
when estimating degrees of freedom for {\it fixed residuals}, our approach matches the one used in
the MASS library of the R programming language~\cite{VenablesRipley}. 
To our knowledge, this approach has not been used for fitting degrees of freedom in general inverse problems, 
and in fact Lange, Little, and Taylor~\cite{Lange1989}, who first proposed Student's t inversion, advocated 
a very different (EM-type) approach for degrees of freedom fitting. 

From a theoretical point of view, the method we propose can be used to solve 
a variety of inverse problems from the general class~\eref{InverseProblemClass}. 
From a practical point of view, the main selling point of the proposed approach is the ability to modify 
existing methods to solve for nuisance parameters on the fly. 

\section{Acknowledgements}
This work was in part financially supported by the Natural Sciences and Engineering Research Council of Canada Discovery Grant (22R81254) 
and the Collaborative Research and Development Grant DNOISE II (375142-08). 
This research was carried out as part of the SINBAD II project with support 
from sponsors of the SINDBAD consortium.
\newpage

\bibliographystyle{plain}
\bibliography{mybib}

\begin{thebibliography}{10}

\bibitem{aravkin12}
A.Y. Aravkin, J.V. Burke, and G.~Pillonetto.
\newblock Robust and trend-following student's t kalman smoothers.
\newblock {\em Optimization Online Preprint}, 2012.

\bibitem{AravkinFriedlanderHerrmannVanleeuwen:2011}
A.Y. Aravkin, M.P. Friedlander, F.J. Herrmann, and T.~van Leeuwen.
\newblock \href{http://www.optimization-online.org/DB_HTML/2011/11/3243.html}
  {Robust inversion, dimensionality reduction, randomized sampling}.
\newblock {\em Submitted to Math. Prog.}, November 2011.

\bibitem{RobustBAArxiv}
A.Y. Aravkin, M.~Styer, Z.~Moratto, A.~Nefian, and M.~Broxton.
\newblock \href{http://arxiv.org/abs/1111.1400} {Student's T Robust Bundle
  Adjustment Algorithm}, November 2011.
\newblock http://arxiv.org/abs/1111.1400.

\bibitem{Bell2008ADo}
Bradley~M. Bell and James~V. Burke.
\newblock Algorithmic differentiation of implicit functions and optimal values.
\newblock In Christian~H. Bischof, H.~Martin B{\"u}cker, Paul~D. Hovland, Uwe
  Naumann, and J.~Utke, editors, {\em Advances in Automatic Differentiation},
  pages 67--77. Springer, 2008.

\bibitem{BellBurkeSchumitzky}
Bradley~M. Bell, James~V. Burke, and Alan Schumitzky.
\newblock A relative weighting method for estimating parameters and variances
  in multiple data sets.
\newblock {\em Computational Statistics \& Data Analysis}, 22(2):119--135, July
  1996.

\bibitem{bunks95}
Carey Bunks, Fatimetou~M. Saleck, S.~Zaleski, and G.~Chavent.
\newblock Multiscale seismic waveform inversion.
\newblock {\em Geophysics}, 60(5):1457--1473, 1995.

\bibitem{fahrmeir2010multivariate}
L.~Fahrmeir and G.~Tutz.
\newblock {\em Multivariate Statistical Modelling Based on Generalized Linear
  Models}.
\newblock Springer Series in Statistics. Springer, 2010.

\bibitem{Fahr1998}
Ludwig Fahrmeir, Rita Kunstler, and Seminar~Fur Statistik.
\newblock Penalized likelihood smoothing in robust state space models.
\newblock {\em Metrika}, 49:173--191, 1998.

\bibitem{Flath:2011}
H.~P. Flath, L.~C. Wilcox, V.~Ak\c{c}elik, J.~Hill, B.~van Bloemen~Waanders,
  and O.~Ghattas.
\newblock Fast algorithms for bayesian uncertainty quantification in
  large-scale linear inverse problems based on low-rank partial hessian
  approximations.
\newblock {\em SIAM J. Sci. Comput.}, 33(1):407--432, February 2011.

\bibitem{GolubPereyra2003}
Gene Golub and Victor Pereyra.
\newblock Separable nonlinear least squares: the variable projection method and
  its applications.
\newblock {\em Inverse Problems}, 19(2):R1, 2003.

\bibitem{GolubPereyra}
{G.H.} Golub and {V.} Pereyra.
\newblock The differentiation of pseudo-inverses and nonlinear least squares
  which variables separate.
\newblock {\em SIAM J. Numer. Anal.}, 10(2):413--432, 1973.

\bibitem{horesh11}
L.~Horesh, E.~Haber, and L.~Tenorio.
\newblock Optimal experimental design for the large-scale nonlinear ill-posed
  problem of impedance imaging.
\newblock In {\em large-scale inverse problems and quantification of
  uncertainty}, Wiley series in computational statistics, pages 273--290.
  Wiley, 2011.

\bibitem{Lange1989}
Kenneth~L. Lange, Roderick J.~A. Little, and Jeremy~M.G. Taylor.
\newblock Robust statistical modeling using the t distribution.
\newblock {\em Journal of the American Statistical Association}, 84:881--896,
  1989.

\bibitem{nelder}
J.~A. Nelder and R.~Mead.
\newblock A simplex method for function minimization.
\newblock {\em Computer Journal}, 7:308--313, 1965.

\bibitem{Osborne2007}
M.~R. Osborne.
\newblock Separable least squares, variable projection, and the
  {G}auss-{N}ewton algorithm.
\newblock {\em Electronic Transactions on Numerical Analysis}, 28(2):1--15,
  2007.

\bibitem{Pratt1998}
R~G Pratt, C~Shin, and Gj~Hicks.
\newblock Gauss-newton and full newton methods in frequency-space seismic
  waveform inversion.
\newblock {\em Geophysical Journal International}, 133(2):341Ã362, 1998.

\bibitem{RTR}
R.T. Rockafellar.
\newblock {\em Convex Analysis}.
\newblock Priceton Landmarks in Mathematics. Princeton University Press, 1970.

\bibitem{Schmidt09optimizingcostly}
Mark Schmidt, Ewout Van~Den Berg, Michael~P. Friedlander, and Kevin Murphy.
\newblock Optimizing costly functions with simple constraints: A limited-memory
  projected quasi-newton algorithm.
\newblock In {\em Proc. of Conf. on Artificial Intelligence and Statistics},
  pages 456--463, 2009.

\bibitem{Tarantola1984}
Albert Tarantola.
\newblock Inversion of seismic reflection data in the acoustic approximation.
\newblock {\em Geophysics}, 49(8):1259--1266, 1984.

\bibitem{VenablesRipley}
W.~N. Venables and B.~D. Ripley.
\newblock {\em Modern Applied Statistics with S}.
\newblock Springer, New York, fourth edition, 2002.
\newblock ISBN 0-387-95457-0.

\end{thebibliography}

\newpage

\end{document}